\documentclass{article}
\usepackage{amsmath}
\usepackage{amssymb}
\usepackage{amsthm}
\usepackage{cite}
\usepackage[arrow,curve,matrix,arc,2cell]{xy}
\UseAllTwocells
\begin{document}
\def\e#1\e{\begin{equation}#1\end{equation}}
\def\ea#1\ea{\begin{align}#1\end{align}}
\def\eq#1{{\rm(\ref{#1})}}
\theoremstyle{plain}
\newtheorem{thm}{Theorem}[section]
\newtheorem{prop}[thm]{Proposition}
\theoremstyle{definition}
\newtheorem{dfn}[thm]{Definition}
\newtheorem{ex}[thm]{Example}
\newtheorem{rem}[thm]{Remark}
\def\dim{\mathop{\rm dim}\nolimits}
\def\supp{\mathop{\rm supp}}
\def\Ker{\mathop{\rm Ker}}
\def\Im{\mathop{\rm Im}}
\def\Hom{\mathop{\rm Hom}\nolimits}
\def\Ho{\mathop{\rm Ho}\nolimits}
\def\Ext{\mathop{\rm Ext}\nolimits}
\def\Aut{\mathop{\rm Aut}}
\def\coh{\mathop{\rm coh}}
\def\qcoh{\mathop{\rm qcoh}}
\def\vect{\mathop{\rm vect}}
\def\id{{\mathop{\rm id}\nolimits}}
\def\Iso{{\rm Iso}}
\def\Spec{\mathop{\rm Spec}\nolimits}
\def\MSpec{\mathop{\rm MSpec}\nolimits}
\def\Top{{\mathop{\bf Top}}}
\def\top{{\kern.05em\rm top}}
\def\tr{{\rm tr}}
\def\nt{{\rm nt}}
\def\CRings{{\mathop{\bf C^{\bs\iy}Rings}}}
\def\CRingsfg{{\mathop{\bf C^{\bs\iy}Rings^{fg}}}}
\def\CRingsfp{{\mathop{\bf C^{\bs\iy}Rings^{fp}}}}
\def\CRingsfa{{\mathop{\bf C^{\bs\iy}Rings^{fa}}}}
\def\CRS{{\mathop{\bf C^{\bs\iy}RS}}}
\def\LCRS{{\mathop{\bf LC^{\bs\iy}RS}}}
\def\CSch{{\mathop{\bf C^{\bs\iy}Sch}}}
\def\CSchlfp{{\mathop{\bf C^{\bs\iy}Sch^{lfp}}}}
\def\CSchlf{{\mathop{\bf C^{\bs\iy}Sch^{lf}}}}
\def\CSta{{\mathop{\bf C^{\bs\iy}Sta}}}
\def\DMCSta{{\mathop{\bf DMC^{\bs\iy}Sta}}}
\def\DMCStalf{{\mathop{\bf DMC^{\bs\iy}Sta^{lf}}}}
\def\DMCStalfp{{\mathop{\bf DMC^{\bs\iy}Sta^{lfp}}}}
\def\Sta{{\mathop{\bf Sta}\nolimits}}
\def\ManSta{{\mathop{\bf ManSta}}}
\def\Man{{\mathop{\bf Man}}}
\def\Orb{{\mathop{\bf Orb}}}
\def\Euc{{\mathop{\bf Euc}}}
\def\Sets{{\mathop{\bf Sets}}}
\def\ul{\underline}
\def\bs{\boldsymbol}
\def\ge{\geqslant}
\def\le{\leqslant\nobreak}
\def\O{{\mathcal O}}
\def\oM{{\mathbin{\smash{\,\,\overline{\!\!\mathcal M\!}\,}}}}
\def\R{{\mathbin{\mathbb R}}}
\def\C{{\mathbin{\mathbb C}}}
\def\CP{{\mathbin{\mathbb{CP}}}}
\def\fC{{\mathbin{\mathfrak C}\kern.05em}}
\def\fD{{\mathbin{\mathfrak D}}}
\def\fE{{\mathbin{\mathfrak E}}}
\def\fF{{\mathbin{\mathfrak F}}}
\def\cC{{\mathbin{\cal C}}}
\def\cD{{\mathbin{\cal D}}}
\def\cE{{\mathbin{\cal E}}}
\def\cF{{\mathbin{\cal F}}}
\def\fCmod{{\mathbin{{\mathfrak C}\text{\rm -mod}}}}
\def\fCmodco{{\mathbin{{\mathfrak C}\text{\rm -mod}^{\rm co}}}}
\def\fCmodfp{{\mathbin{{\mathfrak C}\text{\rm -mod}^{\rm fp}}}}
\def\fDmod{{\mathbin{{\mathfrak D}\text{\rm -mod}}}}
\def\fDmodco{{\mathbin{{\mathfrak D}\text{\rm -mod}^{\rm co}}}}
\def\fDmodfp{{\mathbin{{\mathfrak D}\text{\rm -mod}^{\rm fp}}}}
\def\OXmod{{\mathbin{\O_X\text{\rm -mod}}}}
\def\OYmod{{\mathbin{\O_Y\text{\rm -mod}}}}
\def\OZmod{{\mathbin{\O_Z\text{\rm -mod}}}}
\def\OcXmod{{\mathbin{{\O\kern -0.1em}_{\cal X}\text{\rm -mod}}}}
\def\OcYmod{{\mathbin{\O_{\cal Y}\text{\rm -mod}}}}
\def\OcZmod{{\mathbin{{\O\kern -0.1em}_{\cal Z}\text{\rm -mod}}}}
\def\m{{\mathfrak m}}
\def\ue{{\underline{e}}{}}
\def\uf{{\underline{f\!}\,}{}}
\def\ug{{\underline{g\!}\,}{}}
\def\uh{{\underline{h\!}\,}{}}
\def\ui{{\underline{i\kern -0.07em}\kern 0.07em}{}}
\def\um{{\underline{m\kern -0.1em}\kern 0.1em}{}}
\def\us{{\underline{s\kern -0.15em}\kern 0.15em}{}}
\def\ut{{\underline{t\kern -0.1em}\kern 0.1em}{}}
\def\uu{{\underline{u\kern -0.1em}\kern 0.1em}{}}
\def\uv{{\underline{v\!}\,}{}}
\def\uG{{{\underline{G\!}\,}}{}}
\def\uU{{{\underline{U\kern -0.25em}\kern 0.2em}{}}}
\def\utU{{{\underline{\ti U\kern -0.25em}\kern 0.2em}}{}}
\def\uV{{{\underline{V\kern -0.25em}\kern 0.2em}}{}}
\def\utV{{{\underline{\ti V\kern -0.25em}\kern 0.2em}}{}}
\def\utW{{{\underline{\ti W\!\!}\,\,}}{}}
\def\uW{{{\underline{W\!\!}\,\,}}{}}
\def\uX{{{\underline{X\!}\,}}{}}
\def\uY{{{\underline{Y\!\!}\,\,}}{}}
\def\uZ{{{\underline{Z\!}\,}}{}}
\def\cU{{\cal U}}
\def\cV{{\cal V}}
\def\cW{{\cal W}}
\def\cX{{\cal X}}
\def\cY{{\cal Y}}
\def\cZ{{\cal Z}}
\def\upi{{\underline{\pi\!}\,}}
\def\uid{{\underline{\id\kern -0.1em}\kern 0.1em}}
\def\umu{{\underline{\smash{\mu}}\kern 0.1em}}
\def\al{\alpha}
\def\be{\beta}
\def\ga{\gamma}
\def\de{\delta}
\def\la{\lambda}
\def\th{\theta}
\def\ze{\zeta}
\def\si{\sigma}
\def\om{\omega}
\def\De{\Delta}
\def\Om{\Omega}
\def\Ga{\Gamma}
\def\pd{\partial}
\def\ts{\textstyle}
\def\sm{\setminus}
\def\sh{\sharp}
\def\op{\oplus}
\def\od{\odot}
\def\op{\oplus}
\def\ot{\otimes}
\def\hot{{\mathop{\kern .1em\hat\otimes\kern .1em}\nolimits}}
\def\ov{\overline}
\def\bigop{\bigoplus}
\def\bigot{\bigotimes}
\def\iy{\infty}
\def\es{\emptyset}
\def\ra{\rightarrow}
\def\rra{\rightrightarrows}
\def\Ra{\Rightarrow}
\def\ab{\allowbreak}
\def\longra{\longrightarrow}
\def\hookra{\hookrightarrow}
\def\t{\times}
\def\ci{\circ}
\def\ti{\tilde}
\def\d{{\rm d}}
\def\md#1{\vert #1 \vert}
\def\bmd#1{\big\vert #1 \big\vert}
\def\ms#1{\vert #1 \vert^2}
\def\an#1{\langle #1 \rangle}
\def\ban#1{\bigl\langle #1 \bigr\rangle}
\title{An introduction to $C^\iy$-schemes \\
and $C^\iy$-algebraic geometry}
\author{Dominic Joyce}
\date{}
\maketitle

\section{Introduction}
\label{cs1}

If $X$ is a manifold then the $\R$-algebra $C^\iy(X)$ of smooth
functions $c:X\ra\R$ is a $C^\iy$-{\it ring}. That is, for each
smooth function $f:\R^n\ra\R$ there is an $n$-fold operation
$\Phi_f:C^\iy(X)^n\ra C^\iy(X)$ acting by
$\Phi_f:c_1,\ldots,c_n\mapsto f(c_1,\ldots,c_n)$, and these
operations $\Phi_f$ satisfy many natural identities. Thus,
$C^\iy(X)$ actually has a far richer structure than the obvious
$\R$-algebra structure.

In \cite{Joyc1} the author set out the foundations of a version of
algebraic geometry in which rings or algebras are replaced by
$C^\iy$-rings, focussing on $C^\iy$-{\it schemes}, a category of
geometric objects which generalize manifolds, and whose morphisms
generalize smooth maps, {\it quasicoherent\/} and {\it coherent
sheaves\/} on $C^\iy$-schemes, and $C^\iy$-{\it stacks}, in
particular {\it Deligne--Mumford\/ $C^\iy$-stacks}, a 2-category of
geometric objects which generalize orbifolds. This paper is a survey
of~\cite{Joyc1}.

$C^\iy$-rings and $C^\iy$-schemes were first introduced in {\it
synthetic differential geometry}, see for instance Dubuc
\cite{Dubu2}, Moerdijk and Reyes \cite{MoRe} and Kock \cite{Kock}.
Following Dubuc's discussion of `models of synthetic differential
geometry' \cite{Dubu1} and oversimplifying a bit, symplectic
differential geometers are interested in $C^\iy$-schemes as they
provide a category $\CSch$ of geometric objects which includes
smooth manifolds and certain `infinitesimal' objects, and all fibre
products exist in $\CSch$, and $\CSch$ has some other nice
properties to do with open covers, and exponentials of
infinitesimals.

Synthetic differential geometry concerns proving theorems about
manifolds using synthetic reasoning involving `infinitesimals'. But
one needs to check these methods of synthetic reasoning are valid.
To do this you need a `model', some category of geometric spaces
including manifolds and infinitesimals, in which you can think of
your synthetic arguments as happening. Once you know there exists at
least one model with the properties you want, then as far as
synthetic differential geometry is concerned the job is done. For
this reason $C^\iy$-schemes were not developed very far in synthetic
differential geometry.

Recently, $C^\iy$-rings and $C^\iy$-ringed spaces appeared in a very
different context, as part of Spivak's definition of {\it derived
manifolds\/} \cite{Spiv}, which are an extension to differential
geometry of Jacob Lurie's `derived algebraic geometry' programme.
The author \cite{Joyc2,Joyc3,Joyc4} is developing an alternative
theory of derived differential geometry which simplifies, and goes
beyond, Spivak's derived manifolds. Our notion of derived manifolds
are called {\it d-manifolds}. We also study {\it d-manifolds with
boundary}, and {\it d-manifolds with corners}, and orbifold versions
of all these, {\it d-orbifolds}. To define d-manifolds and
d-orbifolds we need theories of $C^\iy$-schemes, $C^\iy$-stacks, and
quasicoherent sheaves upon them, much of which had not been done, so
the author set up the foundations of these in~\cite{Joyc1}.

D-manifolds and d-orbifolds will have important applications in
symplectic geometry, and elsewhere. Many areas of symplectic
geometry involve moduli spaces $\oM_{g,m}(J,\be)$ of stable
$J$-holomorphic curves in a symplectic manifold $(M,\om)$. The
original motivation for \cite{Joyc2,Joyc3,Joyc4} was to find a good
geometric description for the geometric structure on such moduli
spaces $\oM_{g,m}(J,\be)$. In the Lagrangian Floer cohomology theory
of Fukaya, Oh, Ohta and Ono \cite{FOOO}, moduli spaces
$\oM_{g,m}(J,\be)$ are given the structure of {\it Kuranishi
spaces}. The notion of Kuranishi space seemed to the author to be
unsatisfactory. In trying improve it, using ideas from Spivak
\cite{Spiv}, the author arrived at the theory of
\cite{Joyc2,Joyc3,Joyc4}. The author believes the `correct'
definition of Kuranishi space in the work of Fukaya et al.
\cite{FOOO} should be that a Kuranishi space is a d-orbifold with
corners.

Section \ref{cs2} explains $C^\iy$-rings and their modules,
\S\ref{cs3} introduces $C^\iy$-schemes, and quasicoherent and
coherent sheaves upon them, and \S\ref{cs4} discusses
$C^\iy$-stacks, particularly Deligne--Mumford $C^\iy$-stacks, their
relation to orbifolds, quasicoherent and coherent sheaves on
Deligne--Mumford $C^\iy$-stacks, and orbifold strata of
Deligne--Mumford $C^\iy$-stacks.
\medskip

\noindent{\it Acknowledgements.} I would like to thank Eduardo Dubuc
and Jacob Lurie for helpful conversations.

\section{$C^\iy$-rings}
\label{cs2}

We begin by explaining the basic objects out of which our theories
are built, $C^\iy$-{\it rings}, or {\it smooth rings}, following
\cite[\S 2, \S 3 \& \S 5]{Joyc1}. Everything in
\S\ref{cs21}--\S\ref{cs22} was already known in synthetic
differential geometry, and can be found in Moerdijk and Reyes
\cite[Ch.~I]{MoRe}, Dubuc \cite{Dubu1,Dubu2,Dubu3} or Kock~\cite[\S
III]{Kock}.

\subsection{Two definitions of $C^\iy$-ring}
\label{cs21}

\begin{dfn} A $C^\iy$-{\it ring\/} is a set $\fC$ together with
operations $\Phi_f:\fC^n\ra\fC$ for all $n\ge 0$ and smooth maps
$f:\R^n\ra\R$, where by convention when $n=0$ we define $\fC^0$ to
be the single point $\{\es\}$. These operations must satisfy the
following relations: suppose $m,n\ge 0$, and $f_i:\R^n\ra\R$ for
$i=1,\ldots,m$ and $g:\R^m\ra\R$ are smooth functions. Define a
smooth function $h:\R^n\ra\R$ by
\begin{equation*}
h(x_1,\ldots,x_n)=g\bigl(f_1(x_1,\ldots,x_n),\ldots,f_m(x_1
\ldots,x_n)\bigr),
\end{equation*}
for all $(x_1,\ldots,x_n)\in\R^n$. Then for all
$(c_1,\ldots,c_n)\in\fC^n$ we have
\begin{equation*}
\Phi_h(c_1,\ldots,c_n)=\Phi_g\bigl(\Phi_{f_1}(c_1,\ldots,c_n),
\ldots,\Phi_{f_m}(c_1,\ldots,c_n)\bigr).
\end{equation*}
We also require that for all $1\le j\le n$, defining
$\pi_j:\R^n\ra\R$ by $\pi_j:(x_1,\ldots,x_n)\mapsto x_j$, we have
$\Phi_{\pi_j}(c_1,\ldots,c_n)=c_j$ for
all~$(c_1,\ldots,c_n)\in\fC^n$.

Usually we refer to $\fC$ as the $C^\iy$-ring, leaving the
operations $\Phi_f$ implicit.

A {\it morphism\/} between $C^\iy$-rings $\bigl(\fC,(\Phi_f)_{
f:\R^n\ra\R\,\,C^\iy}\bigr)$, $\bigl({\mathfrak
D},(\Psi_f)_{f:\R^n\ra\R\,\,C^\iy}\bigr)$ is a map
$\phi:\fC\ra{\mathfrak D}$ such that $\Psi_f\bigl(\phi
(c_1),\ldots,\phi(c_n)\bigr)=\phi\ci\Phi_f(c_1,\ldots,c_n)$ for all
smooth $f:\R^n\ra\R$ and $c_1,\ldots,c_n\in\fC$. We will write
$\CRings$ for the category of $C^\iy$-rings.
\label{cs2def1}
\end{dfn}

Here is the motivating example:

\begin{ex} Let $X$ be a manifold, and write $C^\iy(X)$ for the set
of smooth functions $c:X\ra\R$. For $n\ge 0$ and $f:\R^n\ra\R$
smooth, define $\Phi_f:C^\iy(X)^n\ra C^\iy(X)$ by
\e
\bigl(\Phi_f(c_1,\ldots,c_n)\bigr)(x)=f\bigl(c_1(x),\ldots,c_n(x)\bigr),
\label{cs2eq1}
\e
for all $c_1,\ldots,c_n\in C^\iy(X)$ and $x\in X$. It is easy to see
that $C^\iy(X)$ and the operations $\Phi_f$ form a $C^\iy$-ring.

Now let $f:X\ra Y$ be a smooth map of manifolds. Then pullback
$f^*:C^\iy(Y)\ra C^\iy(X)$ mapping $f^*:c\mapsto c\ci f$ is a
morphism of $C^\iy$-rings. Furthermore, every $C^\iy$-ring morphism
$\phi:C^\iy(Y)\ra C^\iy(X)$ is of the form $\phi=f^*$ for a unique
smooth map~$f:X\ra Y$.

Write $\CRings^{\rm op}$ for the opposite category of $\CRings$,
with directions of morphisms reversed, and $\Man$ for the category
of manifolds without boundary. Then we have a full and faithful
functor $F_\Man^\CRings:\Man\ra\CRings^{\rm op}$ acting by
$F_\Man^\CRings(X)=C^\iy(X)$ on objects and $F_\Man^\CRings(f)=f^*$
on morphisms. This embeds $\Man$ as a full subcategory
of~$\CRings^{\rm op}$.
\label{cs2ex1}
\end{ex}

Note that $C^\iy$-rings are far more general than those coming from
manifolds. For example, if $X$ is any topological space we could
define a $C^\iy$-ring $C^0(X)$ to be the set of {\it continuous\/}
$c:X\ra\R$, with operations $\Phi_f$ defined as in \eq{cs2eq1}. For
$X$ a manifold with $\dim X>0$, the $C^\iy$-rings $C^\iy(X)$ and
$C^0(X)$ are different.

There is a more succinct definition of $C^\iy$-rings using category
theory:

\begin{dfn} Write $\Euc$ for the full subcategory of $\Man$ spanned
by the Euclidean spaces $\R^n$. That is, the objects of $\Euc$ are
the manifolds $\R^n$ for $n=0,1,2,\ldots$, and the morphisms in
$\Euc$ are smooth maps $f:\R^n\ra\R^m$. Write $\Sets$ for the
category of sets. In both $\Euc$ and $\Sets$ we have notions of
(finite) products of objects (that is, $\R^{n+m}=\R^n\t\R^m$, and
products $S\t T$ of sets $S,T$), and products of morphisms. Define a
({\it category-theoretic\/}) $C^\iy$-{\it ring\/} to be a
product-preserving functor~$F:\Euc\ra\Sets$.
\label{cs2def2}
\end{dfn}

Here is how this relates to Definition \ref{cs2def1}. Suppose
$F:\Euc\ra\Sets$ is a product-preserving functor. Define
$\fC=F(\R)$. Then $\fC$ is an object in $\Sets$, that is, a set.
Suppose $n\ge 0$ and $f:\R^n\ra\R$ is smooth. Then $f$ is a morphism
in $\Euc$, so $F(f):F(\R^n)\ra F(\R)=\fC$ is a morphism in $\Sets$.
Since $F$ preserves products $F(\R^n)=F(\R)\t\cdots\t F(\R)=\fC^n$,
so $F(f)$ maps $\fC^n\ra\fC$. We define $\Phi_f:\fC^n\ra\fC$ by
$\Phi_f=F(f)$. Then $\bigl(\fC,(\Phi_f)_{f:\R^n\ra\R
\,\,C^\iy}\bigr)$ is a $C^\iy$ ring.

As in \cite[Prop.~2.5]{Joyc1}, \cite[p.~21--22]{MoRe} we have:

\begin{prop} In the category $\CRings$ of\/ $C^\iy$-rings, all
small colimits exist, and so in particular pushouts and all finite
colimits exist.
\label{cs2prop1}
\end{prop}

\begin{dfn} Let $\fC$ be a $C^\iy$-ring. Then we may give $\fC$ the
structure of a {\it commutative\/ $\R$-algebra}. Define addition
`$+$' on $\fC$ by $c+c'=\Phi_f(c,c')$ for $c,c'\in\fC$, where
$f:\R^2\ra\R$ is $f(x,y)=x+y$. Define multiplication `$\,\cdot\,$'
on $\fC$ by $c\cdot c'=\Phi_g(c,c')$, where $g:\R^2\ra\R$ is
$f(x,y)=xy$. Define scalar multiplication by $\la\in\R$ by $\la
c=\Phi_{\la'}(c)$, where $\la':\R\ra\R$ is $\la'(x)=\la x$. Define
elements 0 and 1 in $\fC$ by $0=\Phi_{0'}(\es)$ and
$1=\Phi_{1'}(\es)$, where $0':\R^0\ra\R$ and $1':\R^0\ra\R$ are the
maps $0':\es\mapsto 0$ and $1':\es\mapsto 1$. One can then show
using the relations on the $\Phi_f$ that all the axioms of a
commutative $\R$-algebra are satisfied. In Example \ref{cs2ex1},
this yields the obvious $\R$-algebra structure on the smooth
functions~$c:X\ra\R$.

An {\it ideal\/} $I$ in $\fC$ is an ideal $I\subset\fC$ in $\fC$
regarded as a commutative $\R$-algebra. Then we make the quotient
$\fC/I$ into a $C^\iy$-ring as follows. If $f:\R^n\ra\R$ is smooth,
define $\Phi_f^I:(\fC/I)^n\ra\fC/I$ by
\begin{equation*}
\bigl(\Phi_f^I(c_1+I,\ldots,c_n+I)\bigr)(x)=f\bigl(c_1(x),\ldots,
c_n(x)\bigr)+I.
\end{equation*}
To show this is well-defined, we must show it is independent of the
choice of representatives $c_1,\ldots,c_n$ in $\fC$ for
$c_1+I,\ldots,c_n+I$ in $\fC/I$. By Hadamard's Lemma there exist
smooth functions $g_i:\R^{2n}\ra\R$ for $i=1,\ldots,n$ with
\begin{equation*}
f(y_1,\ldots,y_n)-f(x_1,\ldots,x_n)=\ts\sum_{i=1}^n(y_i-x_i)
g_i(x_1,\ldots,x_n,y_1,\ldots,y_n)
\end{equation*}
for all $x_1,\ldots,x_n,y_1,\ldots,y_n\in\R$. If $c_1',\ldots,c_n'$
are alternative choices for $c_1,\ab\ldots,\ab c_n$, so that
$c_i'+I=c_i+I$ for $i=1,\ldots,n$ and $c_i'-c_i\in I$, we have
\begin{align*}
f\bigl(c_1'(x),&\ldots,c_n'(x)\bigr)-f\bigl(c_1(x),\ldots,
c_n(x)\bigr)\\
&=\ts\sum_{i=1}^n(c_i'-c_i)g_i\bigl(c'_1(x),\ldots,c'_n(x),
c_1(x),\ldots,c_n(x)\bigr).
\end{align*}
The second line lies in $I$ as $c_i'-c_i\in I$ and $I$ is an ideal,
so $\Phi_f^I$ is well-defined, and clearly $\bigl(\fC/I,(\Phi_f^I)_{
f:\R^n\ra\R\,\,C^\iy}\bigr)$ is a $C^\iy$-ring.

We will use the notation $(f_a:a\in A)$ to denote the ideal in a
$C^\iy$-ring $\fC$ generated by a collection of elements
$f_a\in\fC$, $a\in A$. That is,
\begin{equation*}
(f_a:a\in A)=\bigl\{\ts\sum_{i=1}^nf_{a_i}\cdot c_i:\text{$n\ge 0$,
$a_1,\ldots,a_n\in A$, $c_1,\ldots,c_n\in\fC$}\bigr\}.
\end{equation*}
\label{cs2def3}
\end{dfn}

\subsection{Special classes of $C^\iy$-ring}
\label{cs22}

We define {\it finitely generated}, {\it finitely presented}, {\it
local}, and {\it fair\/} $C^\iy$-rings.

\begin{dfn} A $C^\iy$-ring $\fC$ is called {\it finitely
generated\/} if there exist $c_1,\ldots,c_n$ in $\fC$ which generate
$\fC$ over all $C^\iy$-operations. That is, for each $c\in\fC$ there
exists smooth $f:\R^n\ra\R$ with $c=\Phi_f(c_1,\ldots,c_n)$. Given
such $\fC,c_1,\ldots,c_n$, define $\phi:C^\iy(\R^n)\ra\fC$ by
$\phi(f)=\Phi_f(c_1,\ldots,c_n)$ for smooth $f:\R^n\ra\R$, where
$C^\iy(\R^n)$ is as in Example \ref{cs2ex1} with $X=\R^n$. Then
$\phi$ is a surjective morphism of $C^\iy$-rings, so $I=\Ker\phi$ is
an ideal in $C^\iy(\R^n)$, and $\fC\cong C^\iy(\R^n)/I$ as a
$C^\iy$-ring. Thus, $\fC$ is finitely generated if and only if
$\fC\cong C^\iy(\R^n)/I$ for some $n\ge 0$ and ideal $I$
in~$C^\iy(\R^n)$.

An ideal $I$ in $C^\iy(\R^n)$ is called {\it finitely generated\/}
if $I=(f_1,\ldots,f_k)$ for some $f_1,\ldots,f_k\in C^\iy(\R^n)$. A
$C^\iy$-ring $\fC$ is called {\it finitely presented\/} if $\fC\cong
C^\iy(\R^n)/I$ for some $n\ge 0$, where $I$ is a finitely generated
ideal in $C^\iy(\R^n)$.

A difference with conventional algebraic geometry is that
$C^\iy(\R^n)$ is not noetherian, so ideals in $C^\iy(\R^n)$ may not
be finitely generated, and $\fC$ finitely generated does not imply
$\fC$ finitely presented.
\label{cs2def4}
\end{dfn}

\begin{dfn} A $C^\iy$-ring $\fC$ is called a $C^\iy$-{\it local
ring\/} if regarded as an $\R$-algebra, as in Definition
\ref{cs2def3}, $\fC$ is a local $\R$-algebra with residue field
$\R$. That is, $\fC$ has a unique maximal ideal $\m_\fC$
with~$\fC/\m_\fC\cong\R$.

If $\fC,\fD$ are $C^\iy$-local rings with maximal ideals
$\m_\fC,\m_\fD$, and $\phi:\fC\ra\fD$ is a morphism of $C^\iy$
rings, then using the fact that $\fC/\m_\fC\cong\R\cong\fD/\m_\fD$
we see that $\phi^{-1}(\m_\fD)=\m_\fC$, that is, $\phi$ is a {\it
local\/} morphism of $C^\iy$-local rings. Thus, there is no
difference between morphisms and local morphisms.
\label{cs2def5}
\end{dfn}

\begin{ex} For $n\ge 0$ and $p\in\R^n$, define $C^\iy_p(\R^n)$ to
be the set of germs of smooth functions $c:\R^n\ra\R$ at $p\in\R^n$.
That is, $C^\iy_p(\R^n)$ is the quotient of the set of pairs $(U,c)$
with $p\in U\subset\R^n$ open and $c:U\ra\R$ smooth by the
equivalence relation $(U,c)\sim(U',c')$ if there exists $p\in
V\subseteq U\cap U'$ open with $c\vert_V\equiv c'\vert_V$. Define
operations $\Phi_f: (C^\iy_p(\R^n))^m\ra C^\iy_p(\R^n)$ for
$f:\R^m\ra\R$ smooth by \eq{cs2eq1}. Then $C^\iy_p(\R^n)$ is a
$C^\iy$-local ring, with maximal ideal $\mathfrak
m=\bigl\{[(U,c)]:c(p)=0\bigr\}$.
\label{cs2ex2}
\end{ex}

\begin{dfn} An ideal $I$ in $C^\iy(\R^n)$ is called {\it fair\/} if
for each $f\in C^\iy(\R^n)$, $f$ lies in $I$ if and only if
$\pi_p(f)$ lies in $\pi_p(I)\subseteq C^\iy_p(\R^n)$ for all
$p\in\R^n$, where $C^\iy_p(\R^n)$ is as in Example \ref{cs2ex2} and
$\pi_p:C^\iy(\R^n)\ra C^\iy_p(\R^n)$ is the natural projection
$\pi_p:c\mapsto[(\R^n,c)]$. A $C^\iy$-ring $\fC$ is called {\it
fair\/} if it is isomorphic to $C^\iy(\R^n)/I$, where $I$ is a fair
ideal.
\label{cs2def6}
\end{dfn}

Our term `fair' was introduced in \cite{Joyc1} for brevity, but the
idea was already well-known. They were introduced by Dubuc
\cite[Def.~11]{Dubu2} under the name `$C^\iy$-rings of finite type
presented by an ideal of local character', and in more recent work
would be called `finitely generated and germ-determined
$C^\iy$-rings'.

As in \cite[\S 2]{Joyc1}, if $C^\iy(\R^m)/I\cong C^\iy(\R^n)/J$ then
$I$ is finitely generated, or fair, if and only if $J$ is. Thus, to
decide whether a $C^\iy$-ring $\fC$ is finitely presented, or fair,
it is enough to test one presentation $\fC\cong C^\iy(\R^n)/I$.
Also, $\fC$ finitely presented implies $\fC$ fair implies $\fC$
finitely generated. Write $\CRingsfp,\ab\CRingsfa$ and $\CRingsfg$
for the full subcategories of finitely presented, fair, and finitely
generated $C^\iy$-rings in $\CRings$, respectively. Then
\begin{equation*}
\CRingsfp\subset\CRingsfa\subset\CRingsfg\subset\CRings.
\end{equation*}

From  \cite[Prop.s 2.23 \& 2.25]{Joyc1} we have:

\begin{prop} The subcategories\/ $\CRingsfg,\CRingsfp$ are closed
under pushouts and all finite colimits in\/ $\CRings,$ but\/
$\CRingsfa$ is not. Nonetheless, pushouts and finite colimits exist
in\/ $\CRingsfa,$ though they may not coincide with pushouts and
finite colimits in~$\CRings$.
\label{cs2prop2}
\end{prop}

Given morphisms $\phi:\fC\ra\fD$ and $\psi:\fC\ra\fE$ in $\CRings$,
the pushout $\fD\amalg_{\phi,\fC,\psi}\fE$ in $\CRings$ should be
thought of as a {\it completed tensor product\/} $\fD\hot_\fC\fE$.
The tensor product $\fD\ot_\fC\fE$ is an $\R$-algebra, but in
general not a $C^\iy$-ring, and to get a $C^\iy$-ring we must take a
completion $\fD\hot_\fC\fE$. When $\fC=\R$, the trivial
$C^\iy$-ring, the pushout $\fD\amalg_\R\fE$ is the coproduct
$\fD\amalg\fE=\fD\hot_\R\fE$. For example, one can show
that~$C^\iy(\R^m)\hot_\R C^\iy(\R^n)\cong C^\iy(\R^{m+n})$.

Here is~\cite[Prop.~3.1]{Joyc1}.

\begin{prop}{\bf(a)} If\/ $X$ is a manifold without boundary then
the $C^\iy$-ring $C^\iy(X)$ of Example\/ {\rm\ref{cs2ex1}} is
finitely presented.
\smallskip

\noindent{\bf(b)} If\/ $X$ is a manifold with boundary, or with
corners, and\/ $\pd X\ne\es,$ then the $C^\iy$-ring $C^\iy(X)$ of
Example\/ {\rm\ref{cs2ex1}} is fair, but is not finitely presented.
\label{cs2prop3}
\end{prop}

To save space we will say no more about manifolds with boundary or
corners and $C^\iy$-geometry in this paper. More information can be
found in~\cite{Joyc1,Joyc2,Joyc3,Joyc4}.

\begin{ex} A {\it Weil algebra\/} \cite[Def.~1.4]{Dubu1} is a
finite-dimensional commutative $\R$-algebra $W$ which has a maximal
ideal $\m$ with $W/\m\cong\R$ and $\m^n=0$ for some $n>0$. Then by
Dubuc \cite[Prop.~1.5]{Dubu1} or Kock \cite[Th.~III.5.3]{Kock},
there is a unique way to make $W$ into a $C^\iy$-ring compatible
with the given underlying commutative $\R$-algebra. This
$C^\iy$-ring is finitely presented \cite[Prop.~III.5.11]{Kock}.
$C^\iy$-rings from Weil algebras are important in synthetic
differential geometry, in arguments involving infinitesimals.
\label{cs2ex3}
\end{ex}

\subsection{Modules over $C^\iy$-rings, and cotangent modules}
\label{cs23}

In \cite[\S 5]{Joyc1} we discuss modules over $C^\iy$-rings.

\begin{dfn} Let $\fC$ be a $C^\iy$-ring. A $\fC$-{\it module} $M$
is a module over $\fC$ regarded as a commutative $\R$-algebra as in
Definition \ref{cs2def3}. $\fC$-modules form an abelian category,
which we write as $\fCmod$. For example, $\fC$ is a $\fC$-module,
and more generally $\fC\ot_\R V$ is a $\fC$-module for any real
vector space~$V$.

A $\fC$-module $M$ is called {\it finitely presented\/} if there
exists an exact sequence $\fC\ot_\R\R^m\ra\fC\ot_\R\R^n\ra M\ra 0$
in $\fCmod$ for some $m,n\ge 0$. We write $\fCmod^{\rm fp}$ for the
full subcategory of finitely presented $\fC$-modules in $\fCmod$.
Then $\fCmodfp$ is closed under cokernels and extensions in
$\fCmod$. But it may not be closed under kernels, so $\fCmodfp$ may
not be an abelian category.

Let $\phi:\fC\ra\fD$ be a morphism of $C^\iy$-rings. If $M$ is a
$\fC$-module then $\phi_*(M)=M\ot_\fC\fD$ is a $\fD$-module. This
induces a functor $\phi_*:\fCmod\ra\fDmod$, which
maps~$\fCmodfp\ra\fDmodfp$.
\label{cs2def7}
\end{dfn}

\begin{ex} Let $X$ be a manifold, and $E\ra X$ be a vector bundle.
Write $C^\iy(E)$ for the vector space of smooth sections $e$ of $E$.
Then $C^\iy(X)$ acts on $C^\iy(E)$ by $(c,e)\mapsto c\cdot e$ for
$c\in C^\iy(X)$ and $e\in C^\iy(E)$, so $C^\iy(E)$ is a
$C^\iy(X)$-module, which is finitely presented.

Now let $X,Y$ be manifolds and $f:X\ra Y$ a smooth map. Then
$f^*:C^\iy(Y)\ra C^\iy(X)$ is a morphism of $C^\iy$-rings. If $E$ is
a vector bundle over $Y$, then $f^*(E)$ is a vector bundle over $X$.
Under the functor $(f^*)_*:C^\iy(Y)$-mod$\,\ra C^\iy(X)$-mod of
Definition \ref{cs2def7}, we see that $(f^*)_*\bigl(C^\iy(E)
\bigr)=C^\iy(E)\ot_{C^\iy(Y)}C^\iy(X)$ is isomorphic as a
$C^\iy(X)$-module to~$C^\iy\bigl(f^*(E)\bigr)$.
\label{cs2ex4}
\end{ex}

Every commutative algebra $A$ has a natural module $\Om_A$ called
the {\it module of K\"ahler differentials}, which is a kind of
analogue for $A$ of the cotangent bundle $T^*X$ of a manifold $X$.
In \cite[\S 5.3]{Joyc1} we define the {\it cotangent module\/}
$\Om_\fC$ of a $C^\iy$-ring $\fC$, which is the $C^\iy$-version of
the module of K\"ahler differentials.

\begin{dfn} Let $\fC$ be a $C^\iy$-ring, and $M$ a $\fC$-module. A
$C^\iy$-{\it derivation} is an $\R$-linear map $\d:\fC\ra M$ such
that whenever $f:\R^n\ra\R$ is a smooth map and
$c_1,\ldots,c_n\in\fC$, we have
\begin{equation*}
\d\Phi_f(c_1,\ldots,c_n)=\ts\sum_{i=1}^n\Phi_{\frac{\pd f}{\pd
x_i}}(c_1,\ldots,c_n)\cdot\d c_i.
\end{equation*}
Note that $\d$ is {\it not\/} a morphism of $\fC$-modules. We call
such a pair $M,\d$ a {\it cotangent module\/} for $\fC$ if it has
the universal property that for any $\fC$-module $M'$ and
$C^\iy$-derivation $\d':\fC\ra M'$, there exists a unique morphism
of $\fC$-modules $\phi:M\ra M'$ with~$\d'=\phi\ci\d$.

Define $\Om_\fC$ to be the quotient of the free $\fC$-module with
basis of symbols $\d c$ for $c\in\fC$ by the $\fC$-submodule spanned
by all expressions of the form
$\d\Phi_f(c_1,\ldots,c_n)-\sum_{i=1}^n \Phi_{\frac{\pd f}{\pd
x_i}}(c_1,\ldots,c_n)\cdot\d c_i$ for $f:\R^n\ra\R$ smooth and
$c_1,\ldots,c_n\in\fC$, and define $\d_\fC:\fC\ra \Om_\fC$ by
$\d_\fC:c\mapsto\d c$. Then $\Om_\fC,\d_\fC$ is a cotangent module
for $\fC$. Thus cotangent modules always exist, and are unique up to
unique isomorphism.

Let $\fC,\fD$ be $C^\iy$-rings with cotangent modules
$\Om_\fC,\d_\fC$, $\Om_\fD,\d_\fD$, and $\phi:\fC\ra\fD$ be a
morphism of $C^\iy$-rings. Then $\phi$ makes $\Om_\fD$ into a
$\fC$-module, and there is a unique morphism
$\Om_\phi:\Om_\fC\ra\Om_\fD$ in $\fC$-mod with
$\d_\fD\ci\phi=\Om_\phi\ci\d_\fC$. This induces a morphism
$(\Om_\phi)_*:\Om_\fC\ot_\fC\fD\ra\Om_\fD$ in $\fD$-mod with
$(\Om_\phi)_*\ci (\d_\fC\ot\id_\fD)=\d_\fD$. If $\phi:\fC\ra\fD$,
$\psi:\fD\ra\fE$ are morphisms of $C^\iy$-rings
then~$\Om_{\psi\ci\phi}=\Om_\psi\ci\Om_\phi$.
\label{cs2def8}
\end{dfn}

\begin{ex} Let $X$ be a manifold. Then the cotangent bundle $T^*X$
is a vector bundle over $X$, so as in Example \ref{cs2ex4} it yields
a $C^\iy(X)$-module $C^\iy(T^*X)$. The exterior derivative
$\d:C^\iy(X)\ra C^\iy(T^*X)$ is a $C^\iy$-derivation. These
$C^\iy(T^*X),\d$ have the universal property in Definition
\ref{cs2def8}, and so form a cotangent module for~$C^\iy(X)$.

Now let $X,Y$ be manifolds, and $f:X\ra Y$ be smooth. Then
$f^*(TY),TX$ are vector bundles over $X$, and the derivative of $f$
is a vector bundle morphism $\d f:TX\ra f^*(TY)$. The dual of this
morphism is $(\d f)^*:f^*(T^*Y)\ra T^*X$. This induces a morphism of
$C^\iy(X)$-modules $((\d f)^*)_*:C^\iy\bigl(f^*(T^*Y)\bigr)\ra
C^\iy(T^*X)$. This $((\d f)^*)_*$ is identified with $(\Om_{f^*})_*$
in Definition \ref{cs2def8} under the natural isomorphism
$C^\iy\bigl(f^*(T^*Y)\bigr)\cong C^\iy(T^*Y)\ot_{C^\iy(Y)}C^\iy(X)$.
\label{cs2ex5}
\end{ex}

Definition \ref{cs2def8} abstracts the notion of cotangent bundle of
a manifold in a way that makes sense for any $C^\iy$-ring. From
\cite[Th.s~5.13 \& 5.16]{Joyc1} we have:

\begin{thm}{\bf(a)} Suppose $\fC$ is a finitely presented\/
$C^\iy$-ring. Then $\Om_\fC$ is a finitely presented\/ $\fC$-module.
\smallskip

\noindent{\bf(b)} Suppose we are given a pushout diagram of finitely
generated\/ $C^\iy$-rings:
\begin{equation*}
\xymatrix@C=60pt@R=11pt{ \fC \ar[r]_\be \ar[d]^\al & \fE \ar[d]_\de \\
\fD \ar[r]^\ga & \fF,}
\end{equation*}
so that\/ $\fF=\fD\amalg_\fC\fE$. Then the following sequence of\/
$\fF$-modules is exact:
\begin{equation*}
\xymatrix@C=18pt{ \Om_\fC\ot_{\mu_\fC,\fC,\ga\ci\al}\fF
\ar[rrr]^(0.5){(\Om_\al)_*\op -(\Om_\be)_*} &&&
{\begin{subarray}{l}\ts \Om_\fD\ot_{\mu_\fD,\fD,\ga}\fF\,\, \op\\
\ts\Om_\fE\ot_{\mu_\fE,\fE,\de}\fF \end{subarray}}
\ar[rrr]^(0.65){(\Om_\ga)_*\op (\Om_\de)_*} &&& \Om_\fF \ar[r] & 0.}
\end{equation*}
Here\/ $(\Om_\al)_*:\Om_\fC\ot_{\mu_\fC,\fC,\ga\ci\al}\fF\ra
\Om_\fD\ot_{\mu_\fD,\fD,\ga}\fF$ is induced by\/ $\Om_\al:\Om_\fC\ra
\Om_\fD,$ and so on.
\label{cs2thm1}
\end{thm}

\section{$C^\iy$-schemes}
\label{cs3}

We now summarize material in \cite[\S 4]{Joyc1} on $C^\iy$-schemes,
and in \cite[\S 6]{Joyc1} on coherent and quasicoherent sheaves on
$C^\iy$-schemes. Much of \S\ref{cs31} goes back to
Dubuc~\cite{Dubu2}.

\subsection{The definition of $C^\iy$-schemes}
\label{cs31}

The basic definitions are modelled on the definitions of schemes in
Hartshorne \cite[\S II.2]{Hart}, but replacing rings by
$C^\iy$-rings throughout.

\begin{dfn} A {\it $C^\iy$-ringed space\/} $\uX=(X,\O_X)$ is a
topological space $X$ with a sheaf $\O_X$ of $C^\iy$-rings on $X$.
That is, for each open set $U\subseteq X$ we are given a $C^\iy$
ring $\O_X(U)$, and for each inclusion of open sets $V\subseteq
U\subseteq X$ we are given a morphism of $C^\iy$-rings
$\rho_{UV}:\O_X(U)\ra\O_X(V)$, called the {\it restriction maps},
and all this data satisfies the usual sheaf axioms~\cite[\S
II.1]{Hart}.

A {\it morphism\/} $\uf=(f,f^\sh):(X,\O_X)\ra(Y,\O_Y)$ of $C^\iy$
ringed spaces is a continuous map $f:X\ra Y$ and a morphism
$f^\sh:f^{-1}(\O_Y)\ra\O_X$ of sheaves of $C^\iy$-rings on $X$.
There is another way to write the data $f^\sh$: since direct image
of sheaves $f_*$ is right adjoint to inverse image $f^{-1}$, there
is a natural bijection
\e
\Hom_X\bigl(f^{-1}(\O_Y),\O_X\bigr)\cong\Hom_Y\bigl(\O_Y,f_*(\O_X)\bigr).
\label{cs3eq1}
\e
Write $f_\sh:\O_Y\ra f_*(\O_X)$ for the morphism of sheaves of
$C^\iy$-rings on $Y$ corresponding to $f^\sh$ under \eq{cs3eq1}, so
that
\e
f^\sh:f^{-1}(\O_Y)\longra\O_X\quad \leftrightsquigarrow\quad
f_\sh:\O_Y\longra f_*(\O_X).
\label{cs3eq2}
\e
Depending on the application, either $f^\sh$ or $f_\sh$ may be more
useful. We choose to regard $f^\sh$ as primary and write morphisms
as $\uf=(f,f^\sh)$ rather than $(f,f_\sh)$, because in \cite{Joyc2}
we find it convenient to work uniformly using pullbacks, rather than
mixing pullbacks and pushforwards.

A {\it local\/ $C^\iy$-ringed space\/} $\uX=(X,\O_X)$ is a
$C^\iy$-ringed space for which the stalks $\O_{X,x}$ of $\O_X$ at
$x$ are $C^\iy$-local rings for all $x\in X$. Since morphisms of
$C^\iy$-local rings are automatically local morphisms, morphisms of
local $C^\iy$-ringed spaces $(X,\O_X),(Y,\O_Y)$ are just morphisms
of $C^\iy$-ringed spaces, without any additional locality condition.
Write $\CRS$ for the category of $C^\iy$-ringed spaces, and $\LCRS$
for the full subcategory of local $C^\iy$-ringed spaces.

For brevity, we will use the notation that underlined upper case
letters $\uX,\uY,\uZ,\ldots$ represent $C^\iy$-ringed spaces
$(X,\O_X),(Y,\O_Y),(Z,\O_Z),\ldots$, and underlined lower case
letters $\uf,\ug,\ldots$ represent morphisms of $C^\iy$-ringed
spaces $(f,f^\sh),\ab(g,g^\sh),\ab\ldots.$ When we write `$x\in\uX$'
we mean that $\uX=(X,\O_X)$ and $x\in X$. When we write `$\uU$ {\it
is open in\/} $\uX$' we mean that $\uU=(U,\O_U)$ and $\uX=(X,\O_X)$
with $U\subseteq X$ an open set and~$\O_U=\O_X\vert_U$.
\label{cs3def1}
\end{dfn}

\begin{dfn} Write $\CRings^{\rm op}$ for the opposite category of
$\CRings$. The {\it global sections functor\/}
$\Ga:\LCRS\!\ra\!\CRings^{\rm op}$ acts on objects $(X,\O_X)$ in
$\LCRS$ by $\Ga:(X,\O_X)\mapsto\O_X(X)$ and on morphisms
$(f,f^\sh):(X,\O_X)\ra(Y,\O_Y)$ by $\Ga:(f,f^\sh)\mapsto f_\sh(Y)$,
for $f_\sh:\O_X\ra f_*(\O_Y)$ as in \eq{cs3eq2}. As in
\cite[Th.~8]{Dubu2} there is a {\it spectrum functor\/}
$\Spec:\CRings^{\rm op}\ra\LCRS$, defined explicitly in
\cite[Def.~4.12]{Joyc1}, which is a right adjoint to $\Ga$, that is,
for all $\fC\in\CRings$ and $\uX\in\LCRS$ there are functorial
isomorphisms
\e
\Hom_\CRings(\fC,\Ga(\uX))\cong \Hom_\LCRS(\uX,\Spec\fC).
\label{cs3eq3}
\e
For any $C^\iy$-ring $\fC$ there is a natural morphism of
$C^\iy$-rings $\Phi_\fC:\fC\ra\Ga(\Spec\fC)$ corresponding to
$\uid_\uX$ in \eq{cs3eq3} with $\uX=\Spec\fC$. By
\cite[Th.~13]{Dubu2}, the restriction of $\Spec$ to
$(\CRingsfa)^{\rm op}$ is full and faithful.

A local $C^\iy$-ringed space $\uX$ is called an {\it affine\/
$C^\iy$-scheme} if it is isomorphic in $\LCRS$ to $\Spec\fC$ for
some $C^\iy$-ring $\fC$. We call $\uX$ a {\it finitely presented},
or {\it fair}, affine $C^\iy$-scheme if $X\cong\Spec\fC$ for $\fC$
that kind of $C^\iy$-ring.

Let $\uX=(X,\O_X)$ be a local $C^\iy$-ringed space. We call $\uX$ a
$C^\iy$-{\it scheme\/} if $X$ can be covered by open sets
$U\subseteq X$ such that $(U,\O_X\vert_U)$ is an affine
$C^\iy$-scheme. We call a $C^\iy$-scheme $\uX$ {\it locally fair},
or {\it locally finitely presented}, if $X$ can be covered by open
$U\subseteq X$ with $(U,\O_X\vert_U)$ a fair, or finitely presented,
affine $C^\iy$-scheme, respectively.

Write $\CSchlf,\ab\CSchlfp,\CSch$ for the full subcategories of
locally fair, and locally finitely presented, and all,
$C^\iy$-schemes in $\LCRS$, respectively. We call a $C^\iy$-scheme
$\uX$ {\it separated}, or {\it paracompact}, if the underlying
topological space $X$ is Hausdorff, or paracompact.
\label{cs3def2}
\end{dfn}

\begin{ex} Let $X$ be a manifold. Define a $C^\iy$-ringed space
$\uX=(X,\O_X)$ to have topological space $X$ and $\O_X(U)=C^\iy(U)$
for each open $U\subseteq X$, where $C^\iy(U)$ is the $C^\iy$-ring
of smooth maps $c:U\ra\R$, and if $V\subseteq U\subseteq X$ are open
define $\rho_{UV}:C^\iy(U)\ra C^\iy(V)$ by $\rho_{UV}:c\mapsto
c\vert_V$. Then $\uX=(X,\O_X)$ is a local $C^\iy$-ringed space. It
is canonically isomorphic to $\Spec C^\iy(X)$, and so is an affine
$C^\iy$-scheme. It is locally finitely presented.

Define a functor $F_\Man^\CSch:\Man\ra\CSchlfp\subset\CSch$ by
$F_\Man^\CSch=\Spec\ci F_\Man^\CRings$. Then $F_\Man^\CSch$ is full
and faithful, and embeds $\Man$ as a full subcategory of $\CSch$.
\label{cs3ex1}
\end{ex}

By \cite[Prop.s 4.10, 4.11, 4.25, 4.32, Cor.s 4.18, 4.21 \&
Th.~4.33]{Joyc1} we have:

\begin{thm}{\bf(a)} All finite limits exist in the category\/ $\CRS$.
\smallskip

\noindent{\bf(b)} The full subcategories\/ $\CSchlfp,\CSchlf,
\CSch,\LCRS$ in\/ $\CRS$ are closed under all finite limits in\/
$\CRS$. Hence, fibre products and all finite limits exist in each of
these subcategories.
\smallskip

\noindent{\bf(c)} If\/ $\fC$ is a finitely generated\/ $C^\iy$-ring
then\/ $\Spec\fC$ is a fair affine\/ $C^\iy$-scheme.

\smallskip

\noindent{\bf(d)} Let\/ $(X,\O_X)$ be a finitely presented, or fair,
affine $C^\iy$-scheme, and\/ $U\subseteq X$ be an open subset.
Then\/ $(U,\O_X\vert_U)$ is also a finitely presented, or fair,
affine $C^\iy$-scheme, respectively. However, this does not hold for
general affine $C^\iy$-schemes.
\smallskip

\noindent{\bf(e)} Let\/ $(X,\O_X)$ be a locally finitely presented,
locally fair, or general, $C^\iy$-scheme, and\/ $U\subseteq X$ be
open. Then\/ $(U,\O_X\vert_U)$ is also a locally finitely presented,
or locally fair, or general, $C^\iy$-scheme, respectively.
\smallskip

\noindent{\bf(f)} The functor\/ $F_\Man^\CSch$ takes transverse
fibre products in $\Man$ to fibre products in\/~$\CSch$.
\label{cs3thm1}
\end{thm}

In \cite[Def.~4.34 \& Prop.~4.35]{Joyc1} we discuss {\it partitions
of unity\/} on $C^\iy$-schemes, building on ideas of
Dubuc~\cite{Dubu3}.

\begin{dfn} Let $\uX=(X,\O_X)$ be a $C^\iy$-scheme. Consider a
formal sum $\sum_{a\in A}c_a$, where $A$ is an indexing set and
$c_a\in\O_X(X)$ for $a\in A$. We say $\sum_{a\in A}c_a$ is a {\it
locally finite sum on\/} $\uX$ if $X$ can be covered by open
$U\subseteq X$ such that for all but finitely many $a\in A$ we have
$\rho_{XU}(c_a)=0$ in~$\O_X(U)$.

By the sheaf axioms for $\O_X$, if $\sum_{a\in A}c_a$ is a locally
finite sum there exists a unique $c\in\O_X(X)$ such that for all
open $U\subseteq X$ such that $\rho_{XU}(c_a)=0$ in $\O_X(U)$ for
all but finitely many $a\in A$, we have $\rho_{XU}(c)=\sum_{a\in
A}\rho_{XU}(c_a)$ in $\O_X(U)$, where the sum makes sense as there
are only finitely many nonzero terms. We call $c$ the {\it limit\/}
of $\sum_{a\in A}c_a$, written~$\sum_{a\in A}c_a=c$.

Let $c\in\O_X(X)$. Suppose $V_i\subseteq X$ is open and
$\rho_{XV_i}(c)=0\in\O_X(V_i)$ for $i\in I$, and let
$V=\bigcup_{i\in I}V_i$. Then $V\subseteq X$ is open, and
$\rho_{XV}(c)=0\in\O_X(V)$ as $\O_X$ is a sheaf. Thus taking the
union of all open $V\subseteq X$ with $\rho_{XV}(c)=0$ gives a
unique maximal open set $V_c\subseteq X$ such that
$\rho_{XV_c}(c)=0\in\O_X(V_c)$. Define the {\it support\/} $\supp c$
of $c$ to be $X\sm V_c$, so that $\supp c$ is closed in $X$. If
$U\subseteq X$ is open, we say that $c$ {\it is supported in\/} $U$
if~$\supp c\subseteq U$.

Let $\{U_a:a\in A\}$ be an open cover of $X$. A {\it partition of
unity on\/ $\uX$ subordinate to\/} $\{U_a:a\in A\}$ is
$\{\eta_a:a\in A\}$ with $\eta_a\in\O_X(X)$ supported on $U_a$ for
$a\in A$, such that $\sum_{a\in A}\eta_a$ is a locally finite sum on
$\uX$ with~$\sum_{a\in A}\eta_a=1$.
\label{cs3def3}
\end{dfn}

\begin{prop} Suppose $\uX$ is a separated, paracompact, locally fair
$C^\iy$-scheme, and\/ $\{\uU_a:a\in A\}$ an open cover of\/ $\uX$.
Then there exists a partition of unity $\{\eta_a:a\in A\}$ on $\uX$
subordinate to\/~$\{\uU_a:a\in A\}$.
\label{cs3prop1}
\end{prop}

Here are some differences between ordinary schemes and
$C^\iy$-schemes:

\begin{rem}{\bf(i)} If $A$ is a ring or algebra, then points of the
corresponding scheme $\Spec A$ are prime ideals in $A$. However, if
$\fC$ is a $C^\iy$-ring then (by definition) points of $\Spec\fC$
are maximal ideals in $\fC$ with residue field $\R$, or
equivalently, $\R$-algebra morphisms $x:\fC\ra\R$. This has the
effect that if $X$ is a manifold then points of $\Spec C^\iy(X)$ are
just points of~$X$.
\smallskip

\noindent{\bf(ii)} In conventional algebraic geometry, affine
schemes are a restrictive class. Central examples such as $\CP^n$
are not affine, and affine schemes are not closed under open
subsets, so that $\C^2$ is affine but $\C^2\sm\{0\}$ is not. In
contrast, affine $C^\iy$-schemes are already general enough for many
purposes. For example:
\begin{itemize}
\setlength{\itemsep}{0pt}
\setlength{\parsep}{0pt}
\item All manifolds are affine $C^\iy$-schemes.
\item Open $C^\iy$-subschemes of fair affine $C^\iy$-schemes are
fair and affine.
\item Separated, second countable, locally fair $C^\iy$-schemes
are affine.
\end{itemize}
Affine $C^\iy$-schemes are always separated (Hausdorff), so we need
general $C^\iy$-schemes to include non-Hausdorff behaviour.
\smallskip

\noindent{\bf(iii)} In conventional algebraic geometry the Zariski
topology is too coarse for many purposes, so one has to introduce
the \'etale topology. In $C^\iy$-algebraic geometry there is no need
for this, as affine $C^\iy$-schemes are Hausdorff.
\smallskip

\noindent{\bf(iv)} Even very basic $C^\iy$-rings such as
$C^\iy(\R^n)$ for $n>0$ are not noetherian as $\R$-algebras. So
$C^\iy$-schemes should be compared to non-noetherian schemes in
conventional algebraic geometry.
\label{cs3rem}
\end{rem}

\subsection{Quasicoherent and coherent sheaves on $C^\iy$-schemes}
\label{cs32}

In \cite[\S 6]{Joyc1} we discuss sheaves of modules on
$C^\iy$-schemes.

\begin{dfn} Let $\uX=(X,\O_X)$ be a $C^\iy$-scheme. An $\O_X$-{\it
module\/} $\cE$ on $\uX$ assigns a module $\cE(U)$ over $\O_X(U)$
for each open set $U\subseteq X$, with $\O_X(U)$-action
$\mu_U:\O_X(U)\t\cE(U)\ra\cE(U)$, and a linear map
$\cE_{UV}:\cE(U)\ra\cE(V)$ for each inclusion of open sets
$V\subseteq U\subseteq X$, such that the following commutes:
\begin{equation*}
\xymatrix@R=11pt@C=80pt{ \O_X(U)\t \cE(U) \ar[d]^{\rho_{UV}\t
\cE_{UV}} \ar[r]_{\mu_U} & \cE(U) \ar[d]_{\cE_{UV}} \\
\O_X(V)\t \cE(V)\ar[r]^{\mu_V} & \cE(V),}
\end{equation*}
and all this data $\cE(U),\cE_{UV}$ satisfies the usual sheaf
axioms~\cite[\S II.1]{Hart} .

A {\it morphism of\/ $\O_X$-modules\/} $\phi:\cE\ra\cF$ assigns a
morphism of $\O_X(U)$-modules $\phi(U):\cE(U)\ra\cF(U)$ for each
open set $U\subseteq X$, such that $\phi(V)\ci\cE_{UV}=
\cF_{UV}\ci\phi(U)$ for each inclusion of open sets $V\subseteq
U\subseteq X$. Then $\O_X$-modules form an {\it abelian category},
which we write as~$\OXmod$.

As in \cite[\S 6.2]{Joyc1}, the spectrum functor $\Spec:\CRings^{\rm
op}\ab\ra\CSch$ has a counterpart for modules: if $\fC$ is a
$C^\iy$-ring and $(X,\O_X)=\Spec\fC$ we can define a functor
$\MSpec:\fCmod\ra\OXmod$. Let $\uX=(X,\O_X)$ be a $C^\iy$-scheme,
and $\cE$ an $\O_X$-module. We call $\cE$ {\it quasicoherent\/} if
$\uX$ can be covered by open $\uU$ with $\uU\cong\Spec\fC$ for some
$C^\iy$-ring $\fC$, and under this identification
$\cE\vert_\uU\cong\MSpec M$ for some $\fC$-module $M$. We call $\cE$
{\it coherent\/} if furthermore we can take these $\fC$-modules $M$
to be finitely presented. We call $\cE$ a {\it vector bundle of
rank\/} $n\ge 0$ if $\uX$ may be covered by open $\uU$ such that
$\cE\vert_\uU\cong\O_U\ot_\R\R^n$. Write $\qcoh(\uX),\coh(\uX)$, and
$\vect(\uX)$ for the full subcategories of quasicoherent sheaves,
coherent sheaves, and vector bundles in $\OXmod$, respectively.
\label{cs3def4}
\end{dfn}

\begin{dfn} Let $\uf=(f,f^\sh):(X,\O_X)\ra(Y,\O_Y)$ be a morphism
of $C^\iy$-schemes, and $\cE$ be an $\O_Y$-module. Define the {\it
pullback\/} $\uf^*(\cE)$ by $\uf^*(\cE)=f^{-1}(\cE)
\ot_{f^{-1}(\O_Y)}\O_X$, where the inverse image sheaf $f^{-1}(\cE)$
is a sheaf of modules over the inverse image sheaf of $C^\iy$-rings
$f^{-1}(\O_Y)$ on $X$, and the tensor product uses the morphism
$f^\sh:f^{-1}(\O_Y)\ra\O_X$.

If $\phi:\cE\ra\cF$ is a morphism in $\OYmod$ we have an induced
morphism $\uf^*(\phi)=f^{-1}(\phi)\ot\id_{\O_X}:\uf^*(\cE)\ra
\uf^*(\cF)$ in $\OXmod$. Then $\uf^*:\OYmod\ra\OXmod$ is a right
exact functor between abelian categories, which restricts to a right
exact functor~$\uf^*:\qcoh(\uY)\ra\qcoh(\uX)$.

Pullbacks $\uf^*(\cE)$ are a kind of fibre product, and may be
characterized by a universal property. So they should be regarded as
being {\it unique up to canonical isomorphism}, rather than unique.
We use the Axiom of Choice to choose $\uf^*(\cE)$ for all $\uf,\cE$,
and so speak of `the' pullback $\uf^*(\cE)$. However, it may not be
possible to make these choices functorial in $\uf$. That is, if
$\uf:\uX\ra\uY$, $\ug:\uY\ra\uZ$ are morphisms and $\cE\in\OZmod$
then $(\ug\ci\uf)^*(\cE)$ and $\uf^*(\ug^*(\cE))$ are canonically
isomorphic in $\OXmod$, but may not be equal. We will write
$I_{\uf,\ug}(\cE):(\ug\ci\uf)^*(\cE)\ra\uf^*(\ug^*(\cE))$ for these
canonical isomorphisms. Then $I_{\uf,\ug}:(\ug\ci\uf)^*\Ra
\uf^*\ci\ug^*$ is a natural isomorphism of functors.

Similarly, when $\uf$ is the identity $\uid_\uX:\uX\ra\uX$ and
$\cE\in\OXmod$ we may not have $\uid^*_\uX(\cE)=\cE$, but there is a
canonical isomorphism $\de_\uX(\cE):\uid^*_\uX(\cE)\ra\cE$, and
$\de_\uX:\uid^*_\uX\Ra\id_\OXmod$ is a natural isomorphism of
functors.
\label{cs3def5}
\end{dfn}

\begin{ex} Let $X$ be a manifold, and $\uX$ the associated
$C^\iy$-scheme from Example \ref{cs3ex1}, so that $\O_X(U)=C^\iy(U)$
for all open $U\subseteq X$. Let $E\ra X$ be a vector bundle. Define
an $\O_X$-module $\cE$ on $\uX$ by $\cE(U)=C^\iy(E\vert_U)$, the
smooth sections of the vector bundle $E\vert_U\ra U$, and for open
$V\subseteq U\subseteq X$ define $\cE_{UV}:\cE(U)\ra\cE(V)$ by
$\cE_{UV}:e_U\mapsto e_U\vert_V$. Then $\cE\in\vect(\uX)$ is a
vector bundle on $\uX$, which we think of as a lift of $E$ from
manifolds to $C^\iy$-schemes.

Suppose $f:X\ra Y$ is a smooth map of manifolds, and $\uf:\uX\ra\uY$
is the corresponding morphism of $C^\iy$-schemes. Let $F\ra Y$ be a
vector bundle over $Y$, so that $f^*(F)\ra X$ is a vector bundle
over $X$. Let $\cF\in\vect(\uY)$ be the vector bundle over $\uY$
lifting $F$. Then $\uf^*(\cF)$ is canonically isomorphic to the
vector bundle over $\uX$ lifting~$f^*(F)$.
\label{cs3ex2}
\end{ex}

The next theorem comes from \cite[Cor.~6.11 \& Prop.~6.12]{Joyc1}.
In part (a), the reason $\coh(\uX)$ is not closed under kernels is
that the $C^\iy$-rings we are interested in are generally {\it not
noetherian\/} as commutative $\R$-algebras, and this causes problems
with coherence; in conventional algebraic geometry, one usually only
considers coherent sheaves over noetherian schemes.

\begin{thm} {\bf(a)} Let\/ $\uX$ be a $C^\iy$-scheme. Then\/
$\qcoh(\uX)$ is closed under kernels, cokernels and extensions in
$\OXmod,$ so it is an abelian category. Also $\coh(\uX)$ is closed
under cokernels and extensions in\/ $\OXmod,$ but may not be closed
under kernels in\/ $\OXmod,$ so\/ $\coh(\uX)$ may not be an abelian
category.
\smallskip

\noindent{\bf(b)} Suppose\/ $\uf:\uX\ra\uY$ is a morphism of\/
$C^\iy$-schemes. Then pullback\/ $\uf^*:\OYmod\ra\OXmod$ maps
$\qcoh(\uY)\ra\qcoh(\uX)$ and\/ $\coh(\uY)\ra\coh(\uX)$ and\/
$\vect(\uY)\ra\vect(\uX)$. Also $\uf^*:\qcoh(\uY)\ra\qcoh(\uX)$ is a
right exact functor.
\smallskip

\noindent{\bf(c)} Let\/ $\uX$ be a locally fair\/ $C^\iy$-scheme.
Then every\/ $\O_X$-module\/ $\cE$ on\/ $\uX$ is quasicoherent, that
is,~$\qcoh(\uX)=\OXmod$.
\label{cs3thm2}
\end{thm}

Let $\uX$ be a separated, paracompact, locally fair $C^\iy$-scheme.
Then partitions of unity exist on $\uX$ subordinate to any open
cover by Proposition \ref{cs3prop1}. As in \cite[\S 6.3]{Joyc1},
this shows that quasicoherent sheaves $\cE$ on $\uX$ are {\it fine},
which implies that their cohomology groups $H^i(\cE)$ are zero for
all $i>0$. In \cite[Prop.~6.13]{Joyc1} we deduce an exactness
property for sections of quasicoherent sheaves on~$\uX$:

\begin{prop} Suppose\/ $\uX=(X,\O_X)$ is a separated, paracompact,
locally fair $C^\iy$-scheme, and\/ $\cdots \ra\cE^i\,{\buildrel
\phi^i\over\longra}\,\cE^{i+1}\,{\buildrel\phi^{i+1}\over\longra}
\,\cE^{i+2}\ra\cdots$ an exact sequence in $\qcoh(\uX)$. Then
$\cdots\ra\cE^i(U)\,{\buildrel \phi^i(U)\over\longra}\,
\cE^{i+1}(U)\,{\buildrel\phi^{i+1}(U)\over\longra}
\,\cE^{i+2}(U)\ra\cdots$ is an exact sequence of\/ $\O_X(U)$-modules
for each open\/~$U\subseteq X$.
\label{cs3prop2}
\end{prop}

We define cotangent sheaves, the sheaf version of cotangent modules
in~\S\ref{cs23}.

\begin{dfn} Let $\uX$ be a $C^\iy$-scheme. Define ${\cal P}T^*\uX$
to associate to each open $U\subseteq X$ the cotangent module
$\Om_{\O_X(U)}$, and to each inclusion of open sets $V\subseteq
U\subseteq X$ the morphism of $\O_X(U)$-modules
$\Om_{\rho_{UV}}:\Om_{\O_X(U)}\ra\Om_{\O_X(V)}$ associated to the
morphism of $C^\iy$-rings $\rho_{UV}:\O_X(U)\ra\O_X(V)$. Then ${\cal
P}T^*\uX$ is a presheaf of $\O_X$-modules on $\uX$. Define the {\it
cotangent sheaf\/ $T^*\uX$ of\/} $\uX$ to be the sheafification of
${\cal P}T^*\uX$, as an $\O_X$-module.

Let $\uf:\uX\ra\uY$ be a morphism of $C^\iy$-schemes. Then by
Definition \ref{cs3def5}, $\uf^*\bigl(T^*\uY\bigr)=f^{-1}(T^*\uY)
\ot_{f^{-1}(\O_Y)}\O_X,$ where $T^*\uY$ is the sheafification of the
presheaf $V\mapsto\Om_{\O_Y(V)}$, and $f^{-1}(T^*\uY)$ the
sheafification of the presheaf $U\mapsto\lim_{V\supseteq f(U)}
(T^*\uY)(V)$, and $f^{-1}(\O_Y)$ the sheafification of the presheaf
$U\mapsto\lim_{V\supseteq f(U)}\O_Y(V)$. The three sheafifications
combine into one, so that $\uf^*\bigl(T^*\uY\bigr)$ is the
sheafification of the presheaf ${\cal P}(\uf^*(T^*\uY))$ acting by
\begin{equation*}
U\longmapsto{\cal P}(\uf^*(T^*\uY))(U)=
\ts\lim_{V\supseteq f(U)}\Om_{\O_Y(V)}\ot_{\O_Y(V)}\O_X(U).
\end{equation*}

Define a morphism of presheaves ${\cal P}\Om_\uf:{\cal
P}(\uf^*(T^*\uY))\ra{\cal P}T^*\uX$ on $X$ by
\begin{equation*}
({\cal P}\Om_\uf)(U)=\ts\lim_{V\supseteq f(U)}
(\Om_{\rho_{f^{-1}(V)\,U}\ci f_\sh(V)})_*,
\end{equation*}
where $(\Om_{\rho_{f^{-1}(V)\,U}\ci f_\sh(V)})_*:\Om_{\O_Y(V)}
\ot_{\O_Y(V)}\O_X(U)\ra\Om_{\O_X(U)}=({\cal P} T^*\uX)(U)$ is
constructed as in Definition \ref{cs2def8} from the $C^\iy$-ring
morphisms $f_\sh(V):\O_Y(V)\ra\O_X(f^{-1}(V))$ from $f_\sh:\O_Y\ra
f_*(\O_X)$ corresponding to $f^\sh$ in $\uf$ as in \eq{cs3eq2}, and
$\rho_{f^{-1}(V)\,U}:\O_X(f^{-1}(V))\ra\O_X(U)$ in $\O_X$. Define
$\Om_\uf:\uf^*\bigl(T^*\uY\bigr)\ra T^*\uX$ to be the induced
morphism of the associated sheaves.
\label{cs3def6}
\end{dfn}

\begin{ex} Let $X$ be a manifold, and $\uX$ the associated
$C^\iy$-scheme. Then $T^*\uX$ is a vector bundle on $\uX$, and is
canonically isomorphic to the lift to $C^\iy$-schemes from Example
\ref{cs3ex2} of the cotangent vector bundle $T^*X$ of~$X$.
\label{cs3ex3}
\end{ex}

Here \cite[Th.~6.17]{Joyc1} are some properties of cotangent
sheaves.

\begin{thm}{\bf(a)} Let\/ $\uf:\uX\ra\uY$ and\/ $\ug:\uY\ra\uZ$ be
morphisms of\/ $C^\iy$-schemes. Then
\begin{equation*}
\Om_{\ug\ci\uf}=\Om_\uf\ci \uf^*(\Om_\ug)\ci I_{\uf,\ug}(T^*\uZ)
\end{equation*}
as morphisms $(\ug\ci\uf)^*(T^*\uZ)\ra T^*\uX$ in $\OXmod$. Here
$\Om_\ug:\ug^*(T^*\uZ)\ra T^*\uY$ in $\OYmod,$ so applying $\uf^*$
gives $\uf^*(\Om_\ug):\uf^*(\ug^*(T^*\uZ))\ra \uf^*(T^*\uY)$ in\/
$\OXmod,$ and\/ $I_{\uf,\ug}(T^*\uZ):(\ug\ci\uf)^*(T^*\uZ)
\ra\uf^*(\ug^*(T^*\uZ))$ is as in Definition\/~{\rm\ref{cs3def5}}.
\smallskip

\noindent{\bf(b)} Suppose\/ $\uW,\uX,\uY,\uZ$ are locally fair\/
$C^\iy$-schemes with a Cartesian square
\begin{equation*}
\xymatrix@C=65pt@R=11pt{ \uW \ar[r]_\uf \ar[d]^\ue & \uY \ar[d]_\uh \\
\uX \ar[r]^\ug & \uZ}
\end{equation*}
in $\CSchlf,$ so that\/ $\uW=\uX\t_\uZ\uY$. Then the following is
exact in $\qcoh(\uW)\!:$
\begin{equation*}
\xymatrix@C=16pt{ (\ug\ci\ue)^*(T^*\uZ)
\ar[rrrr]^(0.45){\begin{subarray}{l}\ue^*(\Om_\ug)\ci
I_{\ue,\ug}(T^*\uZ)\op\\ -\uf^*(\Om_\uh)\ci
I_{\uf,\uh}(T^*\uZ)\end{subarray}} &&&&
\ue^*(T^*\uX)\!\op\!\uf^*(T^*\uY) \ar[rr]^(0.63){\Om_\ue\op\Om_\uf}
&& T^*\uW \ar[r] & 0.}
\end{equation*}
\label{cs3thm3}
\end{thm}

\section{$C^\iy$-stacks}
\label{cs4}

In \cite[\S 7--\S 11]{Joyc1} we discuss $C^\iy$-{\it stacks}, which
are related to $C^\iy$-schemes in the same way that Artin stacks and
Deligne--Mumford stacks in algebraic geometry are related to
schemes. Stacks are a rather technical subject which take a lot of
work and many pages to set up properly, so to keep this section
short we will give less detail than in \S\ref{cs2} and~\S\ref{cs3}.

We are most interested in a subclass of $C^\iy$-stacks called {\it
Deligne--Mumford\/ $C^\iy$-stacks}. Here are some of their important
properties:
\begin{itemize}
\setlength{\itemsep}{0pt}
\setlength{\parsep}{0pt}
\item Deligne--Mumford\/ $C^\iy$-stacks are geometric
objects locally modelled on quotients $\uU/G$, for $\uU$ an
affine $C^\iy$-scheme and $G$ a finite group.
\item Deligne--Mumford $C^\iy$-stacks are related to $C^\iy$-schemes
in the same way that orbifolds are related to manifolds.
\item Any $C^\iy$-scheme yields an example of a Deligne--Mumford
$C^\iy$-stack.
\item Deligne--Mumford $C^\iy$-stacks form a 2-category
$\DMCSta$. That is, we have objects $\cX,\cY$, 1-morphisms
$f,g:\cX\ra\cY$, and 2-morphisms $\eta:f\Ra g$. All 2-morphisms
are invertible, that is, they are 2-isomorphisms.

The geometric meaning of 1- and 2-morphisms is not obvious; to
get a feel for it, it helps to consider the case when $\cX,\cY$
are quotients $[\uX/G],[\uY/H]$ for $C^\iy$-schemes $\uX,\uY$
and finite groups $G,H$ acting on $\uX,\uY$. Oversimplifying
somewhat, a 1-morphism $f:[\uX/G]\ra[\uY/H]$ is roughly a pair
$(\uf,\rho)$ where $\rho:G\ra H$ is a group morphism and
$\uf:\uX\ra\uY$ is a morphism of $C^\iy$-schemes with $\uf\ci
\ga=\rho(\ga)\ci\uf$ for all $\ga\in G$. If $f=(\uf,\rho)$ and
$g=(\ug,\si)$ are two such 1-morphisms, then a 2-morphism
$\eta:f\Ra g$ is roughly an element $\de\in H$ such that
$\si(\ga)=\de\,\rho(\ga)\de^{-1}$ for all $\ga\in G$,
and~$\ug=\de\ci\uf$.
\item There is a good notion of {\it fibre product\/} in a
2-category. All fibre products of Deligne--Mumford
$C^\iy$-stacks exist, as Deligne--Mumford $C^\iy$-stacks.
\end{itemize}

\subsection{The definition of $C^\iy$-stacks}
\label{cs41}

The next few definitions assume a lot of standard material from
stack theory, which is summarized in~\cite[\S 7]{Joyc1}.

\begin{dfn} Define a {\it Grothendieck topology\/} $\cal J$ on the
category $\CSch$ of $C^\iy$-schemes to have coverings
$\{\ui_a:\uU_a\ra\uU\}_{a\in A}$ where $V_a=i_a(U_a)$ is open in $U$
with $\ui_a:\uU_a\ra(V_a,\O_U\vert_{V_a})$ an isomorphism for all
$a\in A$, and $U=\bigcup_{a\in A}V_a$. Up to isomorphisms of the
$\uU_a$, the coverings $\{\ui_a:\uU_a\ra\uU\}_{a\in A}$ of $\uU$
correspond exactly to open covers $\{V_a:a\in A\}$ of $U$. Then
$(\CSch,{\cal J})$ is a {\it site}.

The {\it stacks on\/} $(\CSch,{\cal J})$ form a 2-category
$\Sta_{(\CSch,{\cal J})}$. The site $(\CSch,{\cal J})$ is {\it
subcanonical}. Thus, if $\uX$ is any $C^\iy$-scheme we have an
associated stack on $(\CSch,{\cal J})$ which we write as $\ul{\bar
X\!}\,$. A $C^\iy$-{\it stack\/} is a stack $\cX$ on $(\CSch,{\cal
J})$ such that the diagonal 1-morphism $\De_\cX:\cX\ra\cX\t\cX$ is
representable, and there exists a surjective 1-morphism
$\Pi:\bar\uU\ra\cX$ called an {\it atlas\/} for some $C^\iy$-scheme
$\uU$. $C^\iy$-stacks form a 2-category $\CSta$. All 2-morphisms in
$\CSta$ are invertible, that is, they are 2-isomorphisms.
\label{cs4def1}
\end{dfn}

\begin{rem} So far as the author knows, \cite{Joyc1} is the first
paper to consider stacks on the site $(\CSch,{\cal J})$. Note that
Behrend and Xu \cite[Def.~2.15]{BeXu} use the term `$C^\iy$-stack'
to mean something different, a geometric stack over the site
$(\Man,{\cal J}_\Man)$ of manifolds without boundary with
Grothendieck topology ${\cal J}_\Man$ given by open covers. These
are called `smooth stacks' by Metzler \cite{Metz}.

Write $\ManSta$ for the 2-category of geometric stacks on
$(\Man,{\cal J}_\Man)$, as in \cite{BeXu,Lerm,Metz}. The functor
$F_\Man^\CSch$ of Example \ref{cs3ex1} embeds the site $(\Man,{\cal
J}_\Man)$ into $(\CSch,{\cal J})$. Thus, restricting from
$(\CSch,{\cal J})$ to $(\Man,{\cal J}_\Man)$ defines a natural
truncation 2-functor $F_\CSta^\ManSta:\CSta\ra\ManSta$.

A $C^\iy$-stack $\cal X$ encodes all morphisms $F:\uU\ra{\cal X}$
for $C^\iy$-schemes $\uU$, whereas its image $F_\CSta^\ManSta({\cal
X})$ remembers only morphisms $F:U\ra{\cal X}$ for manifolds $U$.
Thus the truncation functor $F_\CSta^\ManSta$ {\it loses
information}, as it forgets morphisms from $C^\iy$-schemes which are
not manifolds. This includes any information about {\it
nonreduced\/} $C^\iy$-schemes. For our applications in
\cite{Joyc2,Joyc3,Joyc4} this nonreduced information will be
essential, so we must consider stacks on $(\CSch,{\cal J})$ rather
than on~$(\Man,{\cal J}_\Man)$.
\label{cs4rem}
\end{rem}

\begin{dfn} A {\it groupoid object\/} $(\uU,\uV,\us,\ut,\uu,\ui,\um)$ in
$\CSch$, or simply {\it groupoid\/} in $\CSch$, consists of objects
$\uU,\uV$ in $\CSch$ and morphisms $\us,\ut:\uV\ra\uU$, $\uu:\uU\ra
\uV$, $\ui:\uV\ra\uV$ and $\um:\uV\t_{\us,\uU,\ut}\uV\ra\uV$
satisfying the identities
\begin{gather*}
\us\ci \uu=\ut\ci \uu=\id_{\uU},\;\> \us\ci \ui=\ut,\;\>
\ut\ci\ui=\us,\;\> \us\ci\um=\us\ci\upi_2,\;\> \ut\ci\um=\ut\ci\upi_1,\\
\um\ci(\ui\t\uid_{\uV})=\uu\ci\us,\;\> \um\ci (\uid_{\uV}\t\ui)=\uu\ci\ut,\\
\um\ci(\um\t\uid_{\uV})=\um\ci(\uid_{\uV}\t\um):
\uV\t_{\uU}\uV\t_{\uU}\uV\longra \uV,\\
\um\ci(\uid_{\uV}\t
\uu)=\um\ci(\uu\t\uid_{\uV}):\uV=\uV\t_{\uU}\uU\longra \uV.
\end{gather*}

We write groupoids in $\CSch$ as $\uV\rra\uU$ for short, to
emphasize the morphisms $\us,\ut:\uV\ra\uU$. To any such groupoid we
can associate a {\it groupoid stack\/} $[\uV\rra\uU]$, which is a
$C^\iy$-stack. Conversely, if $\cX$ is a $C^\iy$-stack and
$\Pi:\bar\uU\ra\cX$ is an atlas one can construct a groupoid
$\uV\rra\uU$ in $\CSch$, and $\cX$ is equivalent (in the 2-category
sense) to $[\uV\rra\uU]$. Thus, every $C^\iy$-stack is equivalent to
a groupoid stack.

Suppose $\uU$ is a $C^\iy$-scheme and $G$ is a finite group which
acts on the left on $\uU$ by automorphisms, with action
$\umu:G\t\uU\ra\uU$. Then
\e
\bigl(\uU,G\!\t\!\uU,\upi_{\uU},\umu,1\!\t\!\uid_\uU,
(i\!\ci\!\upi_G)\!\t\!\umu,(m\!\ci\!((\upi_G\!\ci\!\upi_1)\!\t
\!(\upi_G\!\ci\!\upi_2)))\!\t\!(\upi_\uU\!\ci\!\upi_2)\bigr)
\label{cs4eq1}
\e
is a groupoid object in $\CSch$, where $1\in G$ is the identity,
$i:G\ra G$ is the inverse map, $m:G\t G\ra G$ is group
multiplication, and in the final morphism $\upi_1,\upi_2$ are the
projections from $(G\t\uU)\t_{\upi_\uU, \uU,\umu}(G\t\uU)$ to the
first and second factors $G\t\uU$. Write $[\uU/G]$ for the groupoid
stack associated to \eq{cs4eq1}. It is a $C^\iy$-stack, which we
call a {\it quotient stack}.
\label{cs4def2}
\end{dfn}

We define some classes of morphisms of $C^\iy$-schemes.

\begin{dfn} Let $\uf=(f,f^\sh):\uX=(X,\O_X)\ra\uY=(Y,\O_Y)$ be a
morphism in $\CSch$. Then:
\begin{itemize}
\setlength{\itemsep}{0pt}
\setlength{\parsep}{0pt}
\item We call $\uf$ an {\it open embedding\/} if $V=f(X)$ is an open
subset in $Y$ and $(f,f^\sh):(X,\O_X)\ra(V,\O_Y\vert_V)$ is an
isomorphism.
\item We call $\uf$ {\it \'etale\/} if each $x\in X$ has an
open neighbourhood $U$ in $X$ such that $V=f(U)$ is open in $Y$
and $(f\vert_U,f^\sh\vert_U):(U,\O_X\vert_U)\ra (V,\O_Y\vert_V)$
is an isomorphism. That is, $\uf$ is a local isomorphism.
\item We call $\uf$ {\it proper\/} if $f:X\ra Y$ is a proper map of
topological spaces, that is, if $S\subseteq Y$ is compact then
$f^{-1}(S)\subseteq X$ is compact.
\item We call $\uf$ {\it separated\/} if $f:X\ra Y$ is a
separated map of topological spaces, that is, $\De_X=\bigl\{
(x,x):x\in X\bigr\}$ is a closed subset of the topological fibre
product $X\t_{f,Y,f}X=\bigl\{(x,x')\in X\t X:f(x)=f(x')\bigr\}$.
\item We call $\uf$ {\it universally closed\/} if whenever
$\ug:\uW\ra\uY$ is a morphism then $\pi_W:X\t_{f,Y,g}W\ra W$ is
a closed map of topological spaces, that is, it takes closed
sets to closed sets.
\end{itemize}
Each one is invariant under base change and local in the target in
$(\CSch,{\cal J})$. Thus, they are also defined for representable
1-morphisms of $C^\iy$-stacks.
\label{cs4def3}
\end{dfn}

\begin{dfn} Let $\cX$ be a $C^\iy$-stack. We say that $\cX$
is {\it separated\/} if the diagonal 1-morphism
$\De_\cX:\cX\ra\cX\t\cX$ is universally closed. If $\cX=\ul{\bar
X\!}\,$ for some $C^\iy$-scheme $\uX=(X,\O_X)$ then $\cX$ is
separated if and only if $\De_X:X\ra X\t X$ is closed, that is, if
and only if $X$ is Hausdorff, so $\uX$ is separated.
\label{cs4def4}
\end{dfn}

\begin{dfn} Let $\cX$ be a $C^\iy$-stack. A $C^\iy$-{\it
substack\/} $\cY\subset\cX$ is a substack of $\cX$ which is also a
$C^\iy$-stack. It has a natural inclusion 1-morphism
$i_\cY:\cY\hookra\cX$. We call $\cY$ an {\it open\/
$C^\iy$-substack\/} of $\cX$ if $i_\cY$ is a representable open
embedding. An {\it open cover\/} $\{\cY_a:a\in A\}$ of $\cX$ is a
family of open $C^\iy$-substacks $\cY_a\subseteq\cX$ with
$\coprod_{a\in A}i_{\cY_a}:\coprod_{a\in A}\cY_a\ra\cX$ surjective.
\label{cs4def5}
\end{dfn}

{\it Deligne--Mumford stacks\/} in algebraic geometry are locally
modelled on quotient stacks $[X/G]$ for $X$ an affine scheme and $G$
a finite group. This motivates:

\begin{dfn} A {\it Deligne--Mumford\/ $C^\iy$-stack\/} is a
$C^\iy$-stack $\cX$ which admits an open cover $\{\cY_a:a\in A\}$
with each $\cY_a$ equivalent to a quotient stack $[\uU_a/G_a]$ for
$\uU_a$ an affine $C^\iy$-scheme and $G_a$ a finite group. We call
$\cX$ a {\it locally fair}, or {\it locally finitely presented},
Deligne--Mumford $C^\iy$-stack if it has such an open cover with
each $\uU_a$ a fair, or finitely presented, affine $C^\iy$-scheme,
respectively. Write $\DMCStalf,\DMCStalfp$ and $\DMCSta$ for the
full 2-subcategories of locally fair, locally finitely presented,
and all, Deligne--Mumford $C^\iy$-stacks in~$\CSta$.
\label{cs4def6}
\end{dfn}

From \cite[Th.s~8.5, 9.10, 9.16 \& Prop.~9.6]{Joyc1} we have:

\begin{thm}{\bf(a)} All fibre products exist in the\/
$2$-category\/~$\CSta$.
\smallskip

\noindent{\bf(b)} $\DMCSta,\DMCStalf$ and\/ $\DMCStalfp$ are closed
under fibre products in\/~$\CSta$.
\smallskip

\noindent{\bf(c)} $\DMCSta,\kern -.1em\DMCStalf\kern -.1em$ and\/
$\DMCStalfp\kern -.1em$ are closed under taking open\/
$C^\iy$-substacks in\/~$\CSta$.
\smallskip

\noindent{\bf(d)} A\/ $C^\iy$-stack\/ $\cX$ is separated and
Deligne--Mumford if and only if it is equivalent to a groupoid
stack\/ $[\uV\rra\uU]$ where $\uU,\uV$ are separated\/
$C^\iy$-schemes, $\us:\uV\ra\uU$ is \'etale, and
$\us\t\ut:\uV\ra\uU\t\uU$ is universally closed.
\smallskip

\noindent{\bf(e)} A\/ $C^\iy$-stack\/ $\cX$ is separated,
Deligne--Mumford and locally fair (or locally finitely presented) if
and only if it is equivalent to some\/ $[\uV\rra\uU]$ with\/
$\uU,\uV$ separated, locally fair (or locally finitely presented)\/
$C^\iy$-schemes, $\us:\uV\ra\uU$ \'etale, and\/
$\us\t\ut:\uV\ra\uU\t\uU$ proper.
\label{cs4thm1}
\end{thm}

A $C^\iy$-stack $\cX$ has an {\it underlying topological
space}~$\cX_\top$.

\begin{dfn} Let $\cX$ be a $C^\iy$-stack. Write $\ul{*}$ for the
point $\Spec\R$ in $\CSch$, and $\bar{\ul{*}}$ for the associated
point in $\CSta$. Define $\cX_\top$ to be the set of 2-isomorphism
classes $[x]$ of 1-morphisms $x:\bar{\ul{*}}\ra\cX$. If
$\cU\subseteq\cX$ is an open $C^\iy$-substack in $\cX$, then
$\cU_\top\subseteq\cX_\top$. Define ${\cal
T}_{\cX_\top}=\bigl\{\cU_\top:\cU\subseteq\cX$ is an open
$C^\iy$-substack in $\cX\bigr\}$. Then $(\cX_\top,{\cal
T}_{\cX_\top})$ is a topological space, which we call the {\it
underlying topological space\/} of $\cX$, and usually write
as~$\cX_\top$.

If $f:\cX\ra\cY$ is a 1-morphism of $C^\iy$-stacks then there is a
natural continuous map $f_\top:\cX_\top\ra\cY_\top$ defined by
$f_\top([x])=[f\ci x]$. If $f,g:\cX\ra\cY$ are 1-morphisms and
$\eta:f\Ra g$ is a 2-isomorphism then $f_\top=g_\top$. Mapping
$\cX\mapsto\cX_\top$, $f\mapsto f_\top$ and 2-morphisms to
identities defines a 2-functor $F_\CSta^\Top:\CSta\ra\Top$, where
the category of topological spaces $\Top$ is regarded as a
2-category with only identity 2-morphisms.

If $\uX=(X,\O_X)$ is a $C^\iy$-scheme, so that $\ul{\bar X\!}\,$ is
a $C^\iy$-stack, then $\ul{\bar X\!}\,_\top$ is naturally
homeomorphic to $X$, and we will identify $\ul{\bar X\!}\,_\top$
with $X$. If $\uf=(f,f^\sh):\uX=(X,\O_X)\ra\uY=(Y,\O_Y)$ is a
morphism of $C^\iy$-schemes, so that $\ul{\bar f\!}\,:\ul{\bar
X\!}\,\ra\bar\uY$ is a 1-morphism of $C^\iy$-stacks, then $\ul{\bar
f\!}\,_\top:\ul{\bar X\!}\,_\top\ra\bar\uY\!_\top$ is~$f:X\ra Y$.

We call a Deligne--Mumford $C^\iy$-stack $\cX$  {\it paracompact} if
the underlying topological space $\cX_\top$ is paracompact.
\label{cs4def7}
\end{dfn}

\begin{dfn} Let $\cX$ be a $C^\iy$-stack, and $[x]\in\cX_\top$.
Pick a representative $x$ for $[x]$, so that $x:\bar{\ul{*}}\ra\cX$
is a 1-morphism. Let $G$ be the group of 2-morphisms $\eta:x\Ra x$.
There is a natural $C^\iy$-scheme $\uG=(G,\O_G)$ with
$\bar\uG\cong\bar{\ul{*}}\t_{x,\cX,x}\bar{\ul{*}}$, which makes
$\uG$ into a $C^\iy$-{\it group\/} (a group object in $\CSch$, just
as a Lie group is a group object in $\Man$). With $[x]$ fixed, this
$C^\iy$-group $\uG$ is independent of choices up to noncanonical
isomorphism; roughly, $\uG$ is canonical up to conjugation in $\uG$.
We define the {\it isotropy group\/} (or {\it stabilizer group}, or
{\it orbifold group\/}) $\Iso_\cX([x])$ of $[x]$ to be this
$C^\iy$-group $\uG$, regarded as a $C^\iy$-group up to noncanonical
isomorphism.

If $f:{\cal X}\ra{\cal Y}$ is a 1-morphism of $C^\iy$-stacks and
$[x]\in{\cal X}_\top$ with $f_\top([x])=[y]\in{\cal Y}_\top$, for
$y=f\ci x$, then we define $f_*:\Iso_\cX([x])\ra\Iso_\cY([y])$ by
$f_*(\eta)=\id_f*\eta$. Then $f_*$ is a group morphism, and extends
to a $C^\iy$-group morphism. It is independent of choices of
$x\in[x]$, $y\in[y]$ up to conjugation
in~$\Iso_\cX([x]),\Iso_\cY([y])$.
\label{cs4def8}
\end{dfn}

If $\cX$ is a Deligne--Mumford $C^\iy$-stack then $\Iso_\cX([x])$ is
a finite group for all $[x]$ in $\cX_\top$, which is discrete as a
$C^\iy$-group. Here are~\cite[Th.s~9.17(a) \& 9.19]{Joyc1}.

\begin{prop} Let\/ $\cX$ be a Deligne--Mumford $C^\iy$-stack\/ and
$[x]\in\cX_\top,$ so that\/ $\Iso_\cX([x])\cong H$ for some finite
group $H$. Then there exists an open $C^\iy$-substack\/ $\cU$ in
$\cX$ with\/ $[x]\in\cU_\top\subseteq\cX_\top$ and an equivalence\/
$\cU\simeq[\uY/H],$ where $\uY=(Y,\O_Y)$ is an affine $C^\iy$-scheme
with an action of\/ $H,$ and\/ $[x]\in\cU_\top\cong Y/H$ corresponds
to a fixed point\/ $y$ of\/ $H$ in\/~$Y$.
\label{cs4prop1}
\end{prop}

\begin{thm} Suppose $\cX$ is a Deligne--Mumford\/ $C^\iy$-stack
with\/ $\Iso_\cX([x])\cong\{1\}$ for all\/ $[x]\in\cX_\top$. Then\/
$\cX$ is equivalent to\/ $\ul{\bar X\!}\,$ for some\/
$C^\iy$-scheme\/~$\uX$.
\label{cs4thm2}
\end{thm}

In conventional algebraic geometry, a stack with all stabilizer
groups trivial is (equivalent to) an {\it algebraic space}, but may
not be a scheme, so the category of algebraic spaces is larger than
the category of schemes. Here algebraic spaces are spaces which are
locally isomorphic to schemes in the \'etale topology, but not
necessarily locally isomorphic to schemes in the Zariski topology.

In contrast, as Theorem \ref{cs4thm2} shows, in $C^\iy$-algebraic
geometry there is no difference between $C^\iy$-schemes and
$C^\iy$-algebraic spaces. This is because in $C^\iy$-geometry the
Zariski topology is already fine enough, as in Remark
\ref{cs3rem}(iii), so we gain no extra generality by passing to the
\'etale topology.

\subsection{Orbifolds as Deligne--Mumford $C^\iy$-stacks}
\label{cs42}

Example \ref{cs3ex1} defined a functor $F_\Man^\CSch:\Man\ra\CSch$
embedding manifolds as a full subcategory of $C^\iy$-schemes.
Similarly, one might expect to define a (2)-functor
$F_\Orb^\DMCSta:\Orb\ra\DMCSta$ embedding the (2-)category of
orbifolds as a full (2-)subcategory of Deligne--Mumford
$C^\iy$-stacks. In fact, in \cite[\S 9.6]{Joyc1} we took a slightly
different approach: we defined a full 2-subcategory $\Orb$ in
$\DMCSta$, and then showed this is equivalent to other definitions
of the (2-)category of orbifolds. The reason for this is that there
is not one definition of orbifolds, but several, and our new
definition of orbifolds as examples of $C^\iy$-stacks may be as
useful as some of the other definitions.

Orbifolds (without boundary) are spaces locally modelled on $\R^n/G$
for $G$ a finite group acting linearly on $\R^n$, just as manifolds
are spaces locally modelled on $\R^n$. They were introduced by
Satake \cite{Sata}, who called them V-manifolds. Moerdijk
\cite{Moer} defines orbifolds as {\it proper \'etale Lie
groupoids\/} in $\Man$. Both \cite{Moer,Sata} regard orbifolds as an
ordinary category, if the question arises at all. However, for
issues such as fibre products or pullbacks of vector bundles this
category is badly behaved, and it becomes clear that orbifolds
should really be a 2-category, as for stacks in algebraic geometry.

There are two main routes in the literature to defining a 2-category
of orbifolds $\Orb$. The first, as in Pronk \cite{Pron} and Lerman
\cite[\S 3.3]{Lerm}, is to define a 2-category $\bf Gpoid$ of proper
\'etale Lie groupoids, and then to define $\Orb$ as a (weak)
2-category localization of $\bf Gpoid$ at a suitable class of
1-morphisms. The second, as in Behrend and Xu \cite[\S 2]{BeXu},
Lerman \cite[\S 4]{Lerm} and Metzler \cite[\S 3.5]{Metz}, is to
define orbifolds as a class of Deligne--Mumford stacks on the site
$(\Man,{\cal J}_\Man)$ of manifolds with Grothendieck topology
${\cal J}_\Man$ coming from open covers. Our approach is similar to
the second route, but defines orbifolds as a class of
$C^\iy$-stacks, that is, as stacks on the site $(\CSch,{\cal J})$
rather than on~$(\Man,{\cal J}_\Man)$.

\begin{dfn} A $C^\iy$-stack $\cX$ is called an {\it orbifold\/}
if it is equivalent to a groupoid stack $[\uV\rra\uU]$ for some
groupoid $(\uU,\uV,\us,\ut,\uu,\ui,\um)$ in $\CSch$ which is the
image under $F_\Man^\CSch$ of a groupoid $(U,V,s,t,u,i,m)$ in
$\Man$, where $s:V\ra U$ is an \'etale smooth map, and $s\times
t:V\ra U\t U$ is a proper smooth map. That is, $\cX$ is the
$C^\iy$-stack associated to a {\it proper \'etale Lie groupoid\/} in
$\Man$. As a $C^\iy$-stack, every orbifold $\cX$ is a separable,
paracompact, locally finitely presented Deligne--Mumford
$C^\iy$-stack. Write $\Orb$ for the full 2-subcategory of orbifolds
in $\CSta$.
\label{cs4def9}
\end{dfn}

Here is \cite[Th.~9.24 \& Cor.~9.25]{Joyc1}. Since equivalent
(2-)categories are considered to be `the same', the moral of Theorem
\ref{cs4thm3} is that our orbifolds are essentially the same objects
as those considered by other recent authors.

\begin{thm} The $2$-category $\Orb$ of orbifolds defined above is
equivalent to the $2$-categories of orbifolds considered as stacks
on $\Man$ defined in Metzler\/ {\rm\cite[\S 3.4]{Metz}} and Lerman\/
{\rm\cite[\S 4]{Lerm},} and also equivalent as a weak\/ $2$-category
to the weak\/ $2$-categories of orbifolds regarded as proper \'etale
Lie groupoids defined in Pronk\/ {\rm\cite{Pron}} and
Lerman\/~{\rm\cite[\S 3.3]{Lerm}}.

Furthermore, the homotopy category $\Ho(\Orb)$ of\/ $\Orb$ (that is,
the category whose objects are objects in\/ $\Orb,$ and whose
morphisms are $2$-isomorphism classes of\/ $1$-morphisms in\/
$\Orb$) is equivalent to the category of orbifolds regarded as
proper \'etale Lie groupoids defined in Moerdijk\/ {\rm\cite{Moer}}.
Transverse fibre products in\/ $\Orb$ agree with the corresponding
fibre products in\/~$\CSta$.
\label{cs4thm3}
\end{thm}

\subsection{Quasicoherent and coherent sheaves on $C^\iy$-stacks}
\label{cs43}

In \cite[\S 10]{Joyc1} the author studied sheaves on
Deligne--Mumford $C^\iy$-stacks.

\begin{dfn} Let $\cX$ be a Deligne--Mumford $C^\iy$-stack. Define
a category ${\cal C}_\cX$ to have objects pairs $(\uU,u)$ where
$\uU$ is a $C^\iy$-scheme and $u:\bar\uU\ra\cX$ is an \'etale
1-morphism, and morphisms $(\uf,\eta):(\uU,u)\ra(\uV,v)$ where
$\uf:\uU\ra\uV$ is an \'etale morphism of $C^\iy$-schemes, and
$\eta:u\Ra v\ci\ul{\bar f\!}\,$ is a 2-isomorphism. If
$(\uf,\eta):(\uU,u)\ra(\uV,v)$ and $(\ug,\ze):(\uV,v)\ra(\uW,w)$ are
morphisms in ${\cal C}_\cX$ then we define the composition
$(\ug,\ze)\ci(\uf,\eta)$ to be $(\ug\ci\uf,\th):(\uU,u)\ra(\uW,w)$,
where $\th$ is the composition of 2-morphisms across the diagram:
\begin{equation*}
\xymatrix@C=20pt@R=7pt{ \bar\uU \ar[dr]^(0.6){\ul{\bar f\!}\,}
\ar@/^/@<1ex>[drrr]_(0.4)u \ar[dd]_{\overline{\ug\ci\uf}}
\dduppertwocell_{}\omit^{}\omit{<-2.5>^{\id}} & \\
& \bar\uV \ar[rr]^(0.5){v} \ar[dl]^(0.45){\bar\ug} &
\ultwocell_{}\omit^{}\omit{^\eta} & \cX.  \\
\bar\uW \ar@/_/@<-1ex>[urrr]^(0.35)w && {}
\ultwocell_{}\omit^{}\omit{^\ze} }
\end{equation*}

Define an $\O_\cX$-{\it module\/} $\cE$ to assign an $\O_U$-module
$\cE(\uU,u)$ on $\uU=(U,\O_U)$ for all objects $(\uU,u)$ in ${\cal
C}_\cX$, and an isomorphism $\cE_{(\uf,\eta)}:\uf^*(\cE(\uV,v))
\ra\cE(\uU,u)$ for all morphisms $(\uf,\eta):(\uU,u)\ra(\uV,v)$ in
${\cal C}_\cX$, such that for all $(\uf,\eta),(\ug,\ze),
(\ug\ci\uf,\th)$ as above the following diagram of isomorphisms of
$\O_U$-modules commutes:
\e
\begin{gathered}
\xymatrix@C=7pt@R=7pt{ (\ug\ci\uf)^*\bigl(\cE(\uW,w)\bigr)
\ar[rrrrr]_(0.4){\cE_{(\ug\ci\uf,\th)}} \ar[dr]_{I_{\uf,\ug}(\cE
(\uW,w))\,\,\,\,\,\,\,\,\,\,\,\,{}} &&&&& \cE(\uU,u), \\
& \uf^*\bigl(\ug^*(\cE(\uW,w)\bigr)
\ar[rrr]^(0.55){\uf^*(\cE_{(\ug,\ze)})} &&&
\uf^*\bigl(\cE(\uV,v)\bigr) \ar[ur]_{{}\,\,\,\cE_{(\uf,\eta)}} }
\end{gathered}
\label{cs4eq2}
\e
for $I_{\uf,\ug}(\cE)$ as in Definition~\ref{cs3def5}.

A {\it morphism of\/ $\O_\cX$-modules\/} $\phi:\cE\ra\cF$ assigns a
morphism of $\O_U$-modules $\phi(\uU,u):\cE(\uU,u)\ra\cF(\uU,u)$ for
each object $(\uU,u)$ in ${\cal C}_\cX$, such that for all morphisms
$(\uf,\eta):(\uU,u)\ra(\uV,v)$ in ${\cal C}_\cX$ the following
commutes:
\begin{equation*}
\xymatrix@C=60pt@R=11pt{ \uf^*\bigl(\cE(\uV,v)\bigr)
\ar[d]_{\uf^*(\phi(\uV,v))} \ar[r]_{\cE_{(\uf,\eta)}} & \cE(\uU,u)
\ar[d]^{\phi(\uU,u)} \\ \uf^*\bigl(\cF(\uV,v)\bigr)
\ar[r]^{\cF_{(\uf,\eta)}} & \cF(\uU,u). }
\end{equation*}

We call $\cE$ {\it quasicoherent}, or {\it coherent}, or a {\it
vector bundle of rank\/} $n$, if $\cE(\uU,u)$ is quasicoherent, or
coherent, or a vector bundle of rank $n$, respectively, for all
$(\uU,u)\in{\cal C}_\cX$. Write $\OcXmod$ for the category of
$\O_\cX$-modules, and $\qcoh(\cX)$, $\coh(\cX)$, $\vect(\cX)$ for
the full subcategories of quasicoherent sheaves, coherent sheaves,
and vector bundles, respectively.
\label{cs4def10}
\end{dfn}

Here are \cite[Prop.~10.3 \& Ex.~10.4]{Joyc1}.

\begin{prop} Let\/ $\cX$ be a Deligne--Mumford\/ $C^\iy$-stack.
Then\/ $\OcXmod$ is an abelian category, and\/ $\qcoh(\cX)$ is
closed under kernels, cokernels and extensions in $\OcXmod,$ so it
is also an abelian category. Also\/ $\coh(\cX)$ is closed under
cokernels and extensions in\/ $\OcXmod,$ but it may not be closed
under kernels in\/ $\OcXmod,$ so may not be abelian. If\/ $\cX$ is
locally fair then\/~$\qcoh(\cX)\!=\!\OcXmod$.
\label{cs4prop2}
\end{prop}

\begin{ex} Let $\uX$ be a $C^\iy$-scheme. Then $\cX=\ul{\bar X\!}\,$ is
a Deligne--Mumford $C^\iy$-stack. We will define an {\it inclusion
functor\/} ${\cal I}_\uX:\OXmod\ra\OcXmod$. Let $\cE$ be an object
in $\OXmod$. If $(\uU,u)$ is an object in ${\cal C}_\cX$ then
$u:\bar\uU\ra\cX=\ul{\bar X\!}\,$ is 1-isomorphic to
$\bar\uu:\bar\uU\ra\ul{\bar X\!}\,$ for some unique morphism
$\uu:\uU\ra\uX$. Define $\cE'(\uU,u)=\uu^*(\cE)$. If
$(\uf,\eta):(\uU,u)\ra(\uV,v)$ is a morphism in ${\cal C}_\cX$ and
$\uu,\uv$ are associated to $u,v$ as above, so that $\uu=\uv\ci\uf$,
then define
\begin{equation*}
\smash{\cE'_{(\uf,\eta)}= I_{\uf,\uv}(\cE)^{-1}:\uf^*(\cE'(\uV,v))=
\uf^*\bigl(\uv^*(\cE)\bigr)\longra (\uv\ci\uf)^*(\cE)=\cE'(\uU,u).}
\end{equation*}
Then \eq{cs4eq2} commutes for all $(\uf,\eta),(\ug,\ze)$, so $\cE'$
is an $\O_\cX$-module.

If $\phi:\cE\ra\cF$ is a morphism of $\O_X$-modules then we define a
morphism $\phi':\cE'\ra\cF'$ in $\OcXmod$ by
$\phi'(\uU,u)=\uu^*(\phi)$ for $\uu$ associated to $u$ as above.
Then defining ${\cal I}_\uX:\cE\mapsto\cE'$, ${\cal
I}_\uX:\phi\mapsto\phi'$ gives a functor $\OXmod\ra\OcXmod$, which
induces equivalences between the categories
$\OXmod,\qcoh(\uX),\coh(\uX)$ defined in \S\ref{cs32} and
$\OcXmod,\qcoh(\cX),\coh(\cX)$ above.
\label{cs4ex1}
\end{ex}

In \cite[\S 10.2]{Joyc1} we explain how to describe sheaves on a
Deligne--Mumford $C^\iy$-stack $\cX$ in terms of sheaves on $\uU$
for an \'etale atlas $\Pi:\bar\uU\ra\cX$ for $\cX$. Here
are~\cite[Def.~10.5 \& Th.~10.6]{Joyc1}.

\begin{dfn} Let $\cX$ be a Deligne--Mumford $C^\iy$-stack. Then
$\cX$ admits an \'etale atlas $\Pi:\bar\uU\ra\cX$, and as in
Definition \ref{cs4def2} from $\Pi$ we can construct a groupoid
$(\uU,\uV,\us,\ut,\uu,\ui,\um)$ in $\CSch$, with $\us,\ut:\uV\ra\uU$
\'etale, such that $\cX$ is equivalent to the groupoid stack
$[\uV\rra\uU]$. Define a $(\uV\rra\uU)$-{\it module\/} to be a pair
$(E,\Phi)$ where $E$ is an $\O_U$-module and
$\Phi:\us^*(E)\ra\ut^*(E)$ is an isomorphism of $\O_V$-modules, such
that
\e
\begin{split}
I_{\um,\ut}(E)^{-1}\ci\um^*(\Phi)\ci I_{\um,\us}(E)=\,
&\bigl(I_{\upi_1,\ut}(E)^{-1}\ci\upi_1^*(\Phi)\ci
I_{\upi_1,\us}(E)\bigr)\ci\\
&\bigl(I_{\upi_2,\ut}(E)^{-1}\ci\upi_2^*(\Phi)\ci
I_{\upi_2,\us}(E)\bigr)
\end{split}
\label{cs4eq3}
\e
in morphisms of $\O_W$-modules $(\us\ci\um)^*(E)\ra(\ut\ci
\um)^*(E)$, where $\uW=\uV\t_{\us,\uU,\ut}\uV$ and
$\upi_1,\upi_2:\uW\ra\uV$ are the projections. Define a {\it
morphism of\/ $(\uV\rra\uU)$-modules\/} $\phi:(E,\Phi)\ra(F,\Psi)$
to be a morphism of $\O_U$-modules $\phi:E\ra F$ such that
$\Psi\ci\us^*(\phi)=\ut^*(\phi)\ci\Phi :\us^*(E)\ra\ut^*(F)$. Then
$(\uV\rra\uU)$-modules form an {\it abelian category\/}
$(\uV\rra\uU)$-mod. Write $\qcoh(\uV\rra\uU)$ and $\coh(\uV\rra\uU)$
for the full subcategories of $(E,\Phi)$ in $(\uV\rra\uU)$-mod with
$E$ quasicoherent, or coherent, respectively. Then
$\qcoh(\uV\rra\uU)$ is abelian. Define a functor
$F_\Pi:\OcXmod\ra(\uV\rra\uU)$-mod by
$F_\Pi:\cE\mapsto\bigl(\cE(\uU,\Pi),\cE_{(\ut,\eta)}^{-1}\ci
\cE_{(\us,\id_{\Pi\ci\us})}\bigr)$ and
$F_\Pi:\phi\mapsto\phi(\uU,\Pi)$. As in \cite[\S 10.2]{Joyc1},
$F_\Pi(\cE)$ does satisfy \eq{cs4eq3} and so lies in
$(\uV\rra\uU)$-mod, and it also maps $\qcoh,\coh(\cX)$
to~$\qcoh,\coh(\uV\rra\uU)$.
\label{cs4def11}
\end{dfn}

\begin{thm} The functor\/ $F_\Pi$ above induces equivalences
between $\OcXmod,\ab\qcoh(\cX),\coh(\cX)$ and $(\uV\rra\uU)${\rm
-mod}$,\qcoh(\uV\rra\uU),\coh(\uV\rra\uU),$ respectively.
\label{cs4thm4}
\end{thm}

For example, if $\cX=[\uY/G]$ for $\uY$ a $C^\iy$-scheme acted on by
a finite group $G$, then Theorem \ref{cs4thm4} shows that
$\qcoh(\cX)$ is equivalent to the abelian category $\qcoh^G(\uY)$ of
$G$-equivariant quasicoherent sheaves on~$\uY$.

In \S\ref{cs32}, for a morphism of $C^\iy$-schemes $\uf:\uX\ra\uY$
we defined a right exact {\it pullback functor\/}
$\uf^*:\OYmod\ra\OXmod$. Pullbacks may not be strictly functorial in
$\uf$, that is, we do not have $\uf^*(\ug^*(\cE))=(\ug\ci\uf)^*
(\cE)$ for all $\uf:\uX\ra\uY$, $\ug:\uY\ra\uZ$ and $\cE\in\OZmod$,
but instead we have canonical isomorphisms
$I_{\uf,\ug}(\cE):(\ug\ci\uf)^*(\cE)\ra\uf^*(\ug^*(\cE))$. We now
generalize this to sheaves on Deligne--Mumford $C^\iy$-stacks. We
must interpret pullback for 2-morphisms as well as 1-morphisms.

\begin{dfn} Let $f:\cX\ra\cY$ be a 1-morphism of
Deligne--Mumford $C^\iy$-stacks, and $\cF$ be an $\O_\cY$-module. A
{\it pullback\/} of $\cF$ to $\cX$ is an $\O_\cX$-module $\cE$,
together with the following data: if $\uU,\uV$ are $C^\iy$-schemes
and $u:\bar\uU\ra\cX$ and $v:\bar\uV\ra\cY$ are \'etale 1-morphisms,
then there is a $C^\iy$-scheme $\uW$ and morphisms
$\upi_\uU:\uW\ra\uU$, $\upi_\uV:\uW\ra\uV$ giving a 2-Cartesian
diagram:
\e
\begin{gathered}
\xymatrix@C=80pt@R=11pt{ \bar\uW \ar[r]_(0.25){\bar\upi_\uV}
\ar[d]_{\bar\upi_\uU}
\drtwocell_{}\omit^{}\omit{^\ze} & \bar\uV \ar[d]^v \\
\bar\uU \ar[r]^(0.7){f\ci u} & \cY.}
\end{gathered}
\label{cs4eq4}
\e
Then an isomorphism
$i(\cF,f,u,v,\ze):\ab\upi^*_\uU\bigl(\cE(\uU,u)\bigr)\ra
\upi^*_\uV\bigl(\cF(\uV,v)\bigr)$ of $\O_W$-modules should be given,
which is functorial in $(\uU,u)$ in ${\cal C}_\cX$ and $(\uV,v)$ in
${\cal C}_\cY$ and the 2-isomorphism $\ze$ in \eq{cs4eq4}. We
usually write pullbacks $\cE$ as $f^*(\cF)$. By
\cite[Prop.~10.9]{Joyc1}, pullbacks $f^*(\cF)$ exist, and are unique
up to unique isomorphism. Using the Axiom of Choice, we choose a
pullback $f^*(\cF)$ for all such $f:\cX\ra\cY$ and~$\cF$.

Let $f:\cX\ra\cY$ be a 1-morphism, and $\phi:\cE\ra\cF$ be a
morphism in $\OcYmod$. Then $f^*(\cE),f^*(\cF)\in\OcXmod$. Define
the {\it pullback morphism\/} $f^*(\phi):f^*(\cE)\ra f^*(\cF)$ to be
the unique morphism in $\OcXmod$ such that whenever
$u:\bar\uU\ra\cX$, $v:\bar\uV\ra\cY$, $\uW,\upi_\uU,\upi_\uV$ are as
above, the following diagram of morphisms of $\O_W$-modules
commutes:
\begin{equation*}
\xymatrix@C=80pt@R=11pt{ \upi^*_\uU\bigl(f^*(\cE)(\uU,u)\bigr)
\ar[r]_{i(\cE,f,u,v,\ze)} \ar[d]_{\pi^*_\uU(f^*(\phi)(\uU,u))} &
\upi^*_\uV\bigl(\cE(\uV,v)\bigr) \ar[d]^{\pi_\uV^*(\phi(\uV,v))} \\
\upi^*_\uU\bigl(f^*(\cF)(\uU,u)\bigr) \ar[r]^{i(\cF,f,u,v,\ze)} &
\upi^*_\uV\bigl(\cF(\uV,v)\bigr).}
\end{equation*}
This defines a functor $f^*:\OcYmod\ra\OcXmod$, which also maps
$\qcoh(\cY)\ra\qcoh(\cX)$ and $\coh(\cY)\ra\coh(\cX)$. It is right
exact by~\cite[Prop.~10.12]{Joyc1}.

Let $f:\cX\ra\cY$ and $g:\cY\ra\cZ$ be 1-morphisms of
Deligne--Mumford $C^\iy$-stacks, and $\cE\in\OcZmod$. Then $(g\ci
f)^*(\cE)$ and $f^*(g^*(\cE))$ both lie in $\OcXmod$. One can show
that $f^*(g^*(\cE))$ is a possible pullback of $\cE$ by $g\ci f$.
Thus as in Definition \ref{cs3def5}, we have a canonical isomorphism
$I_{f,g}(\cE): (g\ci f)^*(\cE)\ra f^*(g^*(\cE))$. This defines a
natural isomorphism of functors~$I_{f,g}:(g\ci f)^*\Ra f^*\ci g^*$.

Let $f,g:\cX\ra\cY$ be 1-morphisms of Deligne--Mumford
$C^\iy$-stacks, $\eta:f\Ra g$ a 2-morphism, and $\cE\in\OcYmod$.
Then we have $\O_\cX$-modules $f^*(\cE),g^*(\cE)$. Define
$\eta^*(\cE):f^*(\cE)\ra g^*(\cE)$ to be the unique isomorphism such
that whenever $\uU,\uV,\uW,u,v,\upi_\uU,\upi_\uV$ are as above, so
that we have 2-Cartesian diagrams
\begin{equation*}
\xymatrix@C=30pt@R=12pt{ \bar\uW \ar[rrr]_(0.8){\bar\upi_\uV}
\ar[d]_{\bar\upi_\uU}  &
\drrtwocell_{}\omit^{}\omit{^{\!\!\!\!\!\!\!\!\!\!\!
\!\!\!\!\!\!\!\!\!\!\!\!\!\ze\od(\eta*\id_{u\ci\bar\upi_\uU})
\,\,\,\,\,\,\,\,\,\,\,\,{}}} && \bar\uV \ar[d]^v & \bar\uW
\ar[rrr]_(0.7){\bar\upi_\uV} \ar[d]_{\bar\upi_\uU}
& \drtwocell_{}\omit^{}\omit{^\ze} && \bar\uV \ar[d]^v \\
\bar\uU \ar[rrr]^(0.8){f\ci u} &&& \cY, & \bar\uU
\ar[rrr]^(0.7){g\ci u} &&& \cY,}
\end{equation*}
as in \eq{cs4eq4}, where in $\ze\od(\eta*\id_{u\ci\bar\upi_\uU})$
`$*$' is horizontal and `$\od$' vertical composition of 2-morphisms,
then we have commuting isomorphisms of $\O_W$-modules:
\begin{equation*}
\xymatrix@C=100pt@R=-1pt{ \upi^*_\uU\bigl(f^*(\cE)(\uU,u)\bigr)
\ar[dr]^{i(\cE,f,u,v,\ze\od(\eta*\id_{u\ci\bar\upi_\uU}))}
\ar[dd]_{\upi^*_\uU((\eta^*(\cE))(\uU,u))} \\
& \upi^*_\uV\bigl(\cE(\uV,v)\bigr). \\
\upi^*_\uU\bigl(g^*(\cE)(\uU,u)\bigr) \ar[ur]_{i(\cE,g,u,v,\ze)}}
\end{equation*}
This defines a natural isomorphism $\eta^*:f^*\Ra g^*$.

If $\cX$ is a Deligne--Mumford $C^\iy$-stack with identity
1-morphism $\id_\cX:\cX\ra\cX$ then for each $\cE\in\OcXmod$, $\cE$
is a possible pullback $\id_\cX^*(\cE)$, so we have a canonical
isomorphism $\de_\cX(\cE): \id_\cX^*(\cE)\ra\cE$. These define a
natural isomorphism~$\de_\cX:\id_\cX^*\Ra\id_{\OcXmod}$.
\label{cs4def12}
\end{dfn}

Here is~\cite[Th.~10.11]{Joyc1}:

\begin{thm} Mapping $\cX$ to $\OcXmod$ for objects\/ $\cX$ in
$\DMCSta,$ and mapping $1$-morphisms\/ $f:\cX\ra\cY$ to\/
$f^*:\OcYmod\ra\OcXmod,$ and mapping\/ $2$-morphisms\/ $\eta:f\Ra g$
to\/ $\eta^*:f^*\Ra g^*$ for $1$-morphisms $f,g:\cX\ra\cY,$ and the
natural isomorphisms\/ $I_{f,g}:(g\ci f)^*\Ra f^*\ci g^*$ for all\/
$1$-morphisms\/ $f:\cX\ra\cY$ and\/ $g:\cY\ra\cZ$ in $\DMCSta,$
and\/ $\de_\cX$ for all\/ $\cX\in\DMCSta,$ together make up a
\begin{bfseries}pseudofunctor\end{bfseries} $(\DMCSta)^{\rm op}\ra
\mathop{\bf AbCat},$ where $\mathop{\bf AbCat}$ is the\/
$2$-category of abelian categories. That is, they satisfy the
conditions:
\begin{itemize}
\setlength{\itemsep}{0pt}
\setlength{\parsep}{0pt}
\item[{\bf(a)}] If\/ $f:\cW\ra\cX,$ $g:\cX\ra\cY,$
$h:\cY\ra\cZ$ are $1$-morphisms in $\DMCSta$ and\/
$\cE\in\OcZmod$ then the following diagram commutes
in\/~$\OcXmod:$
\begin{equation*}
\xymatrix@C=70pt@R=12pt{ (h\ci g\ci f)^*(\cE) \ar[r]_{I_{f,h\ci
g}(\cE)} \ar[d]_{I_{g\ci f,h}(\cE)}
& f^*\bigl((h\ci g)^*(\cE)\bigr) \ar[d]^{f^*(I_{g,h}(\cE))} \\
(g\ci f)^*\bigl(h^*(\cE)\bigr) \ar[r]^{I_{f,g}(h^*(\cE))} &
f^*\bigl(g^*(h^*(\cE))\bigr).}
\end{equation*}
\item[{\bf(b)}] If\/ $f:\cX\ra\cY$ is a $1$-morphism in\/
$\DMCSta$ and\/ $\cE\in\OcYmod$ then the following pairs of
morphisms in\/ $\OcXmod$ are inverse:
\begin{equation*}
\xymatrix@C=7pt@R=10pt{ {\begin{subarray}{l}\ts f^*(\cE)=\\
\ts (f\!\ci\!\id_\cX)^*(\cE)\end{subarray}}
\ar@/^/[rrr]^{I_{\id_\cX,f}(\cE)} &&& \id_\cX^*(f^*(\cE)),
\ar@/^/[lll]^{\de_\cX(f^*(\cE))} & {\begin{subarray}{l}\ts
f^*(\cE)=\\ \ts (\id_\cY\!\ci\! f)^*(\cE)\end{subarray}}
\ar@/^/[rrr]^{I_{f,\id_\cY}(\cE)} &&& f^*(\id_\cY^*(\cE)).
\ar@/^/[lll]^{f^*(\de_\cY(\cE))} }
\end{equation*}
Also $(\id_f)^*(\id_\cE)=\id_{f^*(\cE)}:f^*(\cE)\ra f^*(\cE)$.
\item[{\bf(c)}] If\/ $f,g,h:\cX\ra\cY$ are
$1$-morphisms and\/ $\eta:f\Ra g,$ $\ze:g\Ra h$ are
$2$-morphisms in\/ $\DMCSta,$ so that\/ $\ze\od\eta:f\Ra h$ is
the vertical composition, and\/ $\cE\in\OcYmod,$ then
\begin{equation*}
\ze^*(\cF)\ci\eta^*(\cE)=(\ze\od\eta)^*(\cE):f^*(\cE)\ra
h^*(\cE)\quad\text{in\/ $\OcXmod$.}
\end{equation*}
\item[{\bf(d)}] If\/ $f,\ti f:\cX\ra\cY,$ $g,\ti g:\cY\ra\cZ$
are $1$-morphisms and\/ $\eta:f\Ra f',$ $\ze:g\Ra g'$
$2$-morphisms in\/ $\DMCSta,$ so that\/ $\ze*\eta:g\ci f\Ra \ti
g\ci\ti f$ is the horizontal composition, and\/ $\cE\in\OcZmod,$
then the following commutes in\/~$\OcXmod:$
\begin{equation*}
\xymatrix@C=60pt@R=12pt{ (g\ci f)^*(\cE)
\ar[rr]_{(\ze*\eta)^*(\cE)} \ar[d]_{I_{f,g}(\cE)} &&
(\ti g\ci\ti f)^*(\cE) \ar[d]^{I_{\ti f,\ti g}(\cE)} \\
f^*(g^*(\cE)) \ar[r]^{f^*(\ze^*(\cE))} & f^*(\ti g^*(\cE))
\ar[r]^{\eta^*(\ti g^*(\cE))} & \ti f^*(\ti g^*(\cE)).}
\end{equation*}
\end{itemize}
\label{cs4thm5}
\end{thm}

\begin{dfn} Let $\cX$ be a Deligne--Mumford $C^\iy$-stack. Define
an $\O_\cX$-module $T^*\cX$ called the {\it cotangent sheaf\/} of
$\cX$ by $(T^*\cX)(\uU,u)=T^*\uU$ for all objects $(\uU,u)$ in
${\cal C}_\cX$ and $(T^*\cX)_{(\uf,\eta)}=\Om_\uf:\uf^*(T^*\uV)\ra
T^*\uU$ for all morphisms $(\uf,\eta):(\uU,u)\ra(\uV,v)$ in ${\cal
C}_\cX$, where $T^*\uU$ and $\Om_\uf$ are as in~\S\ref{cs32}.

Let $f:\cX\ra\cY$ be a 1-morphism of Deligne--Mumford
$C^\iy$-stacks. Then $f^*(T^*\cY),T^*\cX$ are $\O_\cX$-modules.
Define $\Om_f:f^*(T^*\cY)\ra T^*\cX$ to be the unique morphism
characterized as follows. Let $u:\bar\uU\ra\cX$, $v:\bar\uV\ra\cY$,
$\uW,\upi_\uU,\upi_\uV$ be as in Definition \ref{cs4def12}, with
\eq{cs4eq4} Cartesian. Then the following diagram of morphisms of
$\O_W$-modules commutes:
\begin{equation*}
\xymatrix@C=17pt@R=15pt{ \upi^*_\uU\bigl(f^*(T^*\cY)(\uU,u)\bigr)
\ar[d]_{\pi^*_\uU(\Om_f(\uU,u))} \ar[rrr]_{i(T^*\cY,f,u,v,\ze)}
&&& \upi^*_\uV\bigl((T^*\cY)(\uV,v)\bigr) \ar@{=}[r] &
\upi^*_\uV(T^*\uV) \ar[d]_{\Om_{\upi_\uV}} \\
\upi^*_\uU\bigl((T^*\cX)(\uU,u)\bigr)
\ar[rrr]^{(T^*\cX)_{(\upi_\uU,\id_{u\ci\upi_\uU})}} &&&
(T^*\cX)(\uW,u\ci\upi_\uU) \ar@{=}[r] & T^*\uW.}
\end{equation*}

If $\Pi:\bar\uU\ra\cX$, $(\uU,\uV,\us,\ut,\uu,\ui,\um)$ and the
functor $F_\Pi:\OcXmod\ra(\uV\rra\uU)$-mod are as in Definition
\ref{cs4def11} then by definition
$F_\Pi(T^*\cX)=(T^*\uU,\Om_\ut^{-1}\ci\Om_\us)$, and so we write
$T^*(\uV\rra\uU)=(T^*\uU,\Om_\ut^{-1}\ci\Om_\us)$
in~$(\uV\rra\uU)$-mod.
\label{cs4def13}
\end{dfn}

Here \cite[Prop.~10.14 \& Th.~10.15]{Joyc1} is the analogue of
Theorem~\ref{cs3thm3}.

\begin{thm}{\bf(a)} Suppose\/ $\cX$ is an $n$-orbifold. Then
$T^*\cX$ is a rank\/ $n$ vector bundle on\/~$\cX$.
\smallskip

\noindent{\bf(b)} Let\/ $f:\cX\ra\cY$ and\/ $g:\cY\ra\cZ$ be\/
$1$-morphisms of Deligne--Mumford\/ $C^\iy$-stacks. Then
\begin{equation*}
\Om_{g\ci f}=\Om_f\ci f^*(\Om_g)\ci I_{f,g}(T^*\cZ)
\end{equation*}
as morphisms\/ $(g\ci f)^*(T^*\cZ)\ra T^*\cX$ in\/~$\OcXmod$.
\smallskip

\noindent{\bf(c)} Let\/ $f,g:\cX\ra\cY$ be\/ $1$-morphisms of
Deligne--Mumford\/ $C^\iy$-stacks and\/ $\eta:f\Ra g$ a
$2$-morphism. Then $\Om_f=\Om_g\ci\eta^*(T^*\cY):f^*(T^*\cY)\ra
T^*\cX$.
\smallskip

\noindent{\bf(d)} Suppose\/ $\cW,\cX,\cY,\cZ$ are locally fair
Deligne--Mumford\/ $C^\iy$-stacks with a $2$-Cartesian square
\begin{equation*}
\xymatrix@C=70pt@R=12pt{ \cW \ar[r]_(0.25)f \ar[d]^e
\drtwocell_{}\omit^{}\omit{^\eta}
 & \cY \ar[d]_h \\
\cX \ar[r]^(0.7)g & \cZ}
\end{equation*}
in $\DMCStalf,$ so that\/ $\cW=\cX\t_\cZ\cY$. Then the following is
exact in $\qcoh(\cW):$
\begin{equation*}
\xymatrix@C=16pt{ (g\ci e)^*(T^*\cZ)
\ar[rrrrrr]^(0.51){\begin{subarray}{l}e^*(\Om_g)\ci I_{e,g}(T^*\cZ)\op\\
-f^*(\Om_h)\ci I_{f,h}(T^*\cZ)\ci\eta^*(T^*\cZ)
\end{subarray}} &&&&&&
{\raisebox{18pt}{$\displaystyle \begin{subarray}{l}\ts e^*(T^*\cX)\op\\
\ts f^*(T^*\cY)\end{subarray}$}} \ar[rr]^(0.58){\Om_e\op \Om_f} &&
T^*\cW \ar[r] & 0.}
\end{equation*}
\label{cs4thm6}
\end{thm}

\subsection{Orbifold strata of Deligne--Mumford $C^\iy$-stacks}
\label{cs44}

Let $\cX$ be a Deligne--Mumford $C^\iy$-stack, with topological
space $\cX_\top$. Then each point $[x]\in\cX_\top$ has an orbifold
group $\Iso_\cX([x])$, a finite group defined up to isomorphism. For
each finite group $\Ga$ we write $\ti{\cal X}{}_{\ci,\top}^\Ga=
\bigl\{[x]\in\cX_\top:\Iso_\cX([x])\cong\Ga\bigr\}$. This is a
locally closed subset of $\cX_\top$, coming from a locally closed
$C^\iy$-substack $\ti{\cal X}{}^\Ga_\ci$ of $\cX$ with inclusion
$\ti O{}_\ci^\Ga(\cX):\ti{\cal X}{}^\Ga_\ci\ra\cX$, with
\e
\cX_\top=\coprod\nolimits_{\substack{\text{isomorphism classes}\\
\text{of finite groups $\Ga$}}} \ti{\cal X}{}_{\ci,\top}^\Ga.
\label{cs4eq5}
\e
One can show that for each $\Ga$, the closure $\,\ov{\!\ti{\cal
X}}{}_{\ci,\top}^{\,\Ga}$ of $\ti{\cal X}{}_{\ci,\top}^\Ga$ in
$\cX_\top$ satisfies
\begin{equation*}
\,\ov{\!\ti{\cal X}}{}_{\ci,\top}^{\,\Ga}\subseteq\coprod
\nolimits_{\substack{
\text{isomorphism classes of finite groups $\De$:}\\
\text{$\Ga$ is isomorphic to a subgroup of
$\De$}}}\ti{\cal X}{}_{\ci,\top}^\De.
\end{equation*}
Thus \eq{cs4eq5} is a {\it stratification\/} of $\cX_\top$, and the
$\ti{\cal X}{}^\Ga_\ci$ are called {\it orbifold strata\/} of~$\cX$.

In \cite[\S 11.1]{Joyc1} we define six variations of this idea,
Deligne--Mumford $C^\iy$-stacks written $\cX^\Ga,\ti{\cal
X}{}^\Ga,\hat{\cal X}{}^\Ga$, and open $C^\iy$-substacks
$\cX{}^\Ga_\ci\subseteq\cX^\Ga$, $\ti{\cal
X}{}^\Ga_\ci\subseteq\ti{\cal X}{}^\Ga$, $\hat{\cal
X}{}^\Ga_\ci\subseteq \hat{\cal X}{}^\Ga$. The geometric points and
orbifold groups of $\cX^\Ga,\ldots,\hat{\cal X}{}^\Ga_\ci$ are given
by:
\begin{itemize}
\setlength{\itemsep}{0pt}
\setlength{\parsep}{0pt}
\item[(i)] Points of $\cX^\Ga$ are isomorphism classes
$[x,\rho]$, where $[x]\in\cX_\top$ and
$\rho:\Ga\ra\Iso_\cX([x])$ is an injective morphism, and
$\Iso_{\smash{\cX^\Ga}}([x,\rho])$ is the centralizer of
$\rho(\Ga)$ in $\Iso_\cX([x])$. Points of
$\cX{}^\Ga_\ci\subseteq\cX^\Ga$ are $[x,\rho]$ with $\rho$ an
isomorphism, and $\Iso_{\smash{ \cX{}^\Ga_\ci}}([x,\rho])\cong
C(\Ga)$, the centre of $\Ga$.
\item[(ii)] Points of $\ti{\cal X}{}^\Ga$ are pairs $[x,\De]$, where
$[x]\in\cX_\top$ and $\De\subseteq\Iso_\cX([x])$ is isomorphic
to $\Ga$, and $\Iso_{\smash{\ti{\cal X}{}^\Ga}}([x,\De])$ is the
normalizer of $\De$ in $\Iso_\cX([x])$. Points of $\ti{\cal
X}{}^\Ga_\ci\subseteq\ti{\cal X}{}^\Ga$ are $[x,\De]$ with
$\De=\Iso_\cX([x])$, and~$\Iso_{\smash{\ti{\cal X}{}^\Ga_\ci}}
([x,\De])\cong\Ga$.
\item[(iii)] Points $[x,\De]$ of $\hat{\cal X}{}^\Ga,
\hat{\cal X}{}^\Ga_\ci$
are the same as for $\ti{\cal X}{}^\Ga,\ti{\cal X}{}^\Ga_\ci$,
but with orbifold groups $\Iso_{\smash{\hat{\cal
X}{}^\Ga}}([x,\De])\cong \Iso_{\smash{\ti{\cal
X}{}^\Ga}}([x,\De])/\De$ and~$\Iso_{\smash{\hat{\cal
X}{}^\Ga_\ci}}([x,\De])\cong\{1\}$.
\end{itemize}
There are 1-morphisms $O^\Ga(\cX),\ldots,\hat\Pi{}^\Ga_\ci(\cX)$
forming a strictly commutative diagram, where the columns are
inclusions of open $C^\iy$-substacks:
\e
\begin{gathered}
\xymatrix@C=55pt@R=6pt{ \cX{}^\Ga_\ci
\ar[rr]^{\ti\Pi{}^\Ga_\ci(\cX)} \ar[dr]_(0.3){O{}^\Ga_\ci(\cX)}
\ar[dd]_\subset \ar@(ul,l)[]_(0.7){\Aut(\Ga)} && \ti{\cal
X}{}^\Ga_\ci \ar[r]^{\hat\Pi{}^\Ga_\ci(\cX)} \ar[dl]^(0.3){\ti
O^\Ga_\ci(\cX)} \ar[dd]^\subset  &
*+[r]{\hat{\cal X}{}^\Ga_\ci\simeq
\bar{\hat\uX}{}^\Ga_\ci\!\!\!\!\!\!}
\ar@<.5ex>[dd]^\subset \\ & \cX \\
\cX^\Ga \ar[rr]_{\ti\Pi^\Ga(\cX)} \ar[ur]^(0.3){O^\Ga(\cX)}
\ar@(dl,l)[]^(0.7){\Aut(\Ga)} && \ti{\cal X}{}^\Ga
\ar[r]_{\hat\Pi{}^\Ga(\cX)} \ar[ul]_(0.3){\ti O^\Ga(\cX)} &
*+[r]{\hat{\cal X}{}^\Ga.} }
\end{gathered}
\label{cs4eq6}
\e
Here $O^\Ga(\cX),\ti O{}^\Ga(\cX),\ti\Pi{}^\Ga(\cX)$ are proper and
representable, and $\hat\Pi{}^\Ga(\cX)$ is proper. Also $\Aut(\Ga)$
acts on $\cX^\Ga,\cX{}^\Ga_\ci$, with $\ti{\cal
X}{}^\Ga\simeq[\cX^\Ga/\Aut(\Ga)]$ and~$\ti{\cal
X}{}^\Ga_\ci\simeq[\cX^\Ga_\ci/\Aut(\Ga)]$.

Note that there are in general no natural 1-morphisms from
$\hat{\cal X}{}^\Ga,\hat{\cal X}{}^\Ga_\ci$ to any of $\cX,\cX^\Ga,
\cX^\Ga_\ci,\ti{\cal X}{}^\Ga,\ti{\cal X}{}^\Ga_\ci$. Although
$\ti{\cal X}{}^\Ga_\ci$ or $\hat\uX{}^\Ga_\ci$ correspond most
closely to the usual idea of orbifold stratum, we will find that
$\cX^\Ga$ and $\ti{\cal X}{}^\Ga$ are most useful in applications to
orbifold and d-orbifold (co)bordism in \cite[Ch.~13]{Joyc2}, in
which it is vital that $O^\Ga(\cX):\cX{}^\Ga\ra\cX$ and $\ti
O{}^\Ga(\ti{\cal X}):\ti{\cal X}{}^\Ga\ra\ti{\cal X}$ are proper.

As in \cite[\S 11.2]{Joyc1}, representable 1-morphisms and their
2-morphisms lift to orbifold strata $\cX^\Ga,\ti{\cal
X}{}^\Ga,\hat{\cal X}{}^\Ga$. That is, if $f:\cX\ra\cY$ is
representable then we have natural representable 1-morphisms
$f^\Ga:\cX^\Ga\ra\cY^\Ga$, $\ti f{}^\Ga:\ti{\cal X}{}^\Ga\ra\ti{\cal
X}{}^\Ga$, $\hat f{}^\Ga:\hat{\cal X}{}^\Ga\ra\hat{\cal X}{}^\Ga$,
and if $\eta:f\Ra g$ is a 2-morphism of representable
$f,g:\cX\ra\cY$ then we have natural 2-morphisms $\eta^\Ga:f^\Ga\Ra
g^\Ga$, $\ti\eta{}^\Ga:\ti f{}^\Ga\Ra\ti g{}^\Ga$,
$\hat\eta{}^\Ga:\hat f{}^\Ga\Ra\hat g{}^\Ga$. These behave in a
strongly functorial way for orbifold strata $\cX^\Ga,\ti{\cal
X}{}^\Ga$, so that $(g\ci f)^\Ga=g^\Ga\ci f^\Ga$, and so on, and in
a weakly functorial way for orbifold strata $\hat{\cal X}{}^\Ga$, so
that $(\,\,\widehat{\!\!g\ci f\!\!}\,\,){}^\Ga$ is 2-isomorphic to
$\hat g{}^\Ga\ci \hat f{}^\Ga$, and so on.

We can describe orbifold strata of quotient
$C^\iy$-stacks,~\cite[Th.~11.9]{Joyc1}:

\begin{thm} Let\/ $\uX$ be a separated\/ $C^\iy$-scheme and\/ $G$ a
finite group acting on $\uX$ by isomorphisms, and write\/
$\cX=[\uX/G]$ for the quotient\/ $C^\iy$-stack, which is a
Deligne--Mumford\/ $C^\iy$-stack. Let\/ $\Ga$ be a finite group.
Then there are equivalences of\/ $C^\iy$-stacks
\ea
\cX^\Ga&\simeq\ts\bigl[\bigl(\coprod_{\text{injective group
morphisms $\rho:\Ga\ra G$}}\uX^{\rho(\Ga)}\bigr)/G\bigr],
\label{cs4eq7}\\
\cX{}^\Ga_\ci&\simeq\ts\bigl[\bigl(\coprod_{\text{injective group
morphisms $\rho:\Ga\ra G$}}\uX^{\rho(\Ga)}_\ci\bigr)/G\bigr],
\label{cs4eq8}\\
\ti{\cal X}{}^\Ga&\simeq\ts\bigl[\bigl(\coprod_{\text{subgroups
$\De\subseteq G$: $\De\cong\Ga$}}\uX^\De\bigr)/G\bigr],
\label{cs4eq9}\\
\ti{\cal X}{}^\Ga_\ci&\simeq\ts\bigl[\bigl(\coprod_{\text{subgroups
$\De\subseteq G$: $\De\cong\Ga$}}\uX^\De_\ci\bigr)/G\bigr],
\label{cs4eq10}\\
\hat{\cal X}{}^\Ga&\simeq\coprod_{\text{conjugacy classes $[\De]$ of
subgroups $\De\subseteq G$ with $\De\cong\Ga$} \!\!\!\!\!\!\!\!\!
\!\!\!\!\!\!\!\!\!\!\!\!\!\!\!\!\!\!\!\!\!\!\!\!\!\!\!\!\!\!\!\!
\!\!\!\!\!\!\!\!\!\!\!\!\!\!\!\!\!\!\!\!\!\!\!}
\bigl[\,\uX^\De\big/\bigl(\{g\in G:\De=g\De
g^{-1}\}/\De\bigr)\bigr],
\label{cs4eq11}\\
\hat{\cal X}{}^\Ga_\ci&\simeq\coprod_{\text{conjugacy classes
$[\De]$ of subgroups $\De\subseteq G$ with $\De\cong\Ga$}
\!\!\!\!\!\!\!\!\!
\!\!\!\!\!\!\!\!\!\!\!\!\!\!\!\!\!\!\!\!\!\!\!\!\!\!\!\!\!\!\!\!
\!\!\!\!\!\!\!\!\!\!\!\!\!\!\!\!\!\!\!\!\!\!\!}
\bigl[\,\uX^\De_\ci\big/\bigl(\{g\in G:\De=g\De
g^{-1}\}/\De\bigr)\bigr].
\label{cs4eq12}
\ea
Here for each subgroup\/ $\De\subseteq G,$ we write\/ $\uX^\De$ for
the closed\/ $C^\iy$-subscheme in\/ $\uX$ fixed by $\De$ in\/ $G,$
and\/ $\uX^\De_\ci$ for the open $C^\iy$-subscheme in $\uX^\De$ of
points in $\uX$ whose stabilizer group in $G$ is exactly\/~$\De$.

Under the equivalences {\rm\eq{cs4eq7}--\eq{cs4eq12},} the
$1$-morphisms in \eq{cs4eq6} are identified up to $2$-isomorphism
with\/ $1$-morphisms between quotient\/ $C^\iy$-stacks induced by
natural\/ $C^\iy$-scheme morphisms between
$\coprod_\rho\uX^{\rho(\Ga)},\uX,\ldots.$ For example, the disjoint
union over $[\rho]$ of the inclusion $\uX^{\rho(\Ga)} \hookra\uX$ is
a $G$-equivariant morphism $\coprod_\rho\uX^{\rho(\Ga)}\ra\uX,$
inducing a $1$-morphism $\smash{[\coprod_\rho\uX^{\rho(\Ga)}/G]\ra
[\uX/G]}$. This is identified with\/ $O^\Ga(\cX):\cX^\Ga\ra\cX$ by\/
\eq{cs4eq7}.
\label{ag11thm2}
\end{thm}

In \cite[\S 11.4--\S 11.6]{Joyc1} we study pullbacks of sheaves to
orbifold strata. We find that if $\cE$ lies in $\qcoh(\cX)$ then the
pullback $\cE^\Ga:=O^\Ga(\cX)^*(\cE)$ to $\cX^\Ga$ has a natural
action of $\Ga$, so we have a splitting
$\cE^\Ga=\cE^\Ga_\tr\op\cE^\Ga_\nt$ into subsheaves with trivial
`$\tr$' and nontrivial `$\nt$' $\Ga$-representations. For cotangent
sheaves we have a natural isomorphism $T^*(\cX^\Ga)\cong
(T^*\cX)^\Ga_\tr$. A similar picture holds for sheaves on orbifold
strata~$\ti{\cal X}{}^\Ga,\hat{\cal X}{}^\Ga$.

\medskip

\noindent{\small\sc The Mathematical Institute, 24-29 St. Giles,
Oxford, OX1 3LB, U.K.}

\noindent{\small\sc E-mail: \tt joyce@maths.ox.ac.uk}

\end{document}